\RequirePackage{fix-cm}
\documentclass[smallextended]{svjour3}       

\smartqed  
\usepackage{graphicx}
\usepackage[utf8]{inputenc}
\usepackage[english]{babel}
\usepackage{listings}   
\usepackage{algorithm2e}
\usepackage{color}
\usepackage {amsmath,amsfonts,amssymb}
\usepackage{subfig}
\usepackage{texnames}
\usepackage{url}
\definecolor{hotpink}{rgb}{0.9,0,0.5}
\usepackage[colorlinks,urlcolor=blue,citecolor=hotpink,linkcolor=blue]{hyperref}
\renewcommand{\leq}{\leqslant}
\renewcommand{\geq}{\geqslant}

\journalname{Journal of Scientific Computing}

\begin{document}


\title{New integration methods for perturbed ODEs  based on symplectic implicit Runge-Kutta schemes with application to solar system simulations
}

\titlerunning{New integration methods based on symplectic implicit Runge-Kutta schemes} 

\author{Mikel Anto\~nana   \and Joseba Makazaga \and Ander Murua
}


\institute{Computer Science and Artificial Intelligence Department, UPV/EHU (University of the Basque Country), Donostia, Spain \\
           \email{Mikel.Antonana@ehu.eus}     
           \and    \email{Joseba.Makazaga@ehu.eus}    
            \and    \email{Ander.Murua@ehu.eus} 
            }

\date{Received: date / Accepted: date}
\maketitle

\begin{abstract}
We propose a family of integrators, Flow-Composed Implicit Runge-Kutta (FCIRK) methods,   for perturbations of nonlinear ordinary differential equations, consisting of the composition of flows of the unperturbed part alternated with one step of an implicit Runge-Kutta (IRK) method applied to a transformed system. The resulting integration schemes are symplectic when both the perturbation and the unperturbed part are Hamiltonian and the underlying IRK scheme is symplectic. In addition, they are symmetric in time (resp. have order of accuracy  $r$) if the underlying IRK scheme is time-symmetric (resp. of order $r$). The proposed new methods admit mixed precision implementation that allows us to efficiently reduce the effect of round-off errors. 
We particularly focus on the potential application to long-term solar system simulations, with the equations of motion of the solar system rewritten as a Hamiltonian perturbation of a system of uncoupled Keplerian equations.
We present some preliminary numerical experiments with a simple point mass Newtonian 10-body model of the solar system (with the sun, the eight planets, and Pluto) written in canonical heliocentric coordinates.
\end{abstract}

\section{Introduction}

We are concerned with the efficient and accurate numerical integration of systems of perturbed ordinary differential equations of the form 
\begin{align}
\label{eq:pertODE}
&\frac{du}{dt} = k(u) + g(u), \\ 
\label{eq:iv}
&u(t_0) = u_0 \in \mathbb{R}^D.
\end{align}
where $g(u)$ is considered as a perturbation, and  the  unperturbed system
\begin{equation}
\label{eq:unpert}
\frac{du}{dt} = k(u),
\end{equation}
can be solved exactly for any initial value. That is, we are able to compute the $t$-flow $\varphi^{[k]}_{t}: \mathbb{R}^{D} \to \mathbb{R}^{D}$ for all $t\in \mathbb{R}$. 

 We propose a family of integration methods intended for the case where the unperturbed system (\ref{eq:unpert}) is non linear.  When (\ref{eq:unpert}) is linear, our integrators reduce to Lawson's Generalized Runge-Kutta methods~\cite{Lawson1967}.
We are particularly interested in long term integrations of problems where both the perturbed system (\ref{eq:pertODE}) and the unperturbed one (\ref{eq:unpert}) are autonomous Hamiltonian systems.

As a motivating example, we consider the N-body equations of a planetary system (in particular, the solar system) with several planets orbiting with near-Keplerian motion around a central massive body. In that case, the corresponding system of ODEs can be written in the form (\ref{eq:pertODE}), where the unperturbed part  (\ref{eq:unpert})  is a Hamiltonian system consisting of several uncoupled Kepler problems, and the perturbation accounts for the interaction among the planets.
With an appropriate choice of the state variables (i.e., Jacobi coordinates, or canonical heliocentric coordinates), both the Keplerian part and the perturbation are Hamiltonian.

Strang splitting method, also known as generalized leapfrog method, is a well known second order integration method with favorable geometric properties that can be successfully applied to systems of the form (\ref{eq:pertODE}). However,   the $t$-flow of the perturbation (in addition to the $t$-flow of the unperturbed system) need to be explicitly computed. 
A similar second order integrator that avoids that requirement can be obtained by replacing the $t$-flow of the perturbation by the approximate discrete flow obtained by applying the implicit midpoint rule to the perturbation.
 More precisely, approximations $u_j \approx u(t_j)$ of the solution $u(t)$ of (\ref{eq:pertODE}) at $t=t_j := t_0 + j h$  (with a fixed step-length $h \in \mathbb{R}$ with small enough absolute value) are computed for $j=1,2,3,\ldots$ as
\begin{equation}
\label{eq:mStrang}
u_{j} = \widehat{\psi}_h(u_{j-1}) :=\varphi_{h/2} (\psi_h(\varphi_{h/2}(u_{j-1}))), 
\end{equation}
%
%
where $\varphi_{h}$ denotes the $h$-flow of the unperturbed system (\ref{eq:unpert}), and the map $\psi_{h}:\mathbb{R}^D \to \mathbb{R}^D$ is implicitly defined as follows: 
\begin{equation}
\label{eq:midpoint}
u^* := \psi_{h}(u) \quad \mbox{such that} \quad
u^* = u + h \, g\left( \textstyle \frac{1}{2} (u + u^* )\right).
\end{equation}
The resulting integration method is symplectic provided that
both the perturbed system (\ref{eq:pertODE}) and the unperturbed one (\ref{eq:unpert}) are Hamiltonian.

In order to understand how the accuracy of a given integrator for perturbed systems varies when the perturbation becomes smaller, it is customary to introduce a small parameter $\epsilon>0$ multiplying the perturbation, that is,
\begin{equation}
\label{eq:pertODEeps}
\frac{du}{dt} = k(u) + \epsilon\, g(u).
\end{equation}
In the case of generalized leapfrog method and its above mentioned modifications, the global error for the numerical solution obtained with time-step length $h$ for a prescribed integration interval behaves like $\mathcal{O}(\epsilon \, h^2)$ as $\epsilon \to 0$ and $h \to 0$.

For high precision long-term integrations, methods of higher order of accuracy are desirable. Here we consider high order methods of the form (\ref{eq:mStrang}) with the map $\psi_{h}:\mathbb{R}^D \to \mathbb{R}^D$ defined as the result of the application of one step of a Runge-Kutta scheme to a transformed system (a different one at each integration step).  Such transformed system is obtained by applying a time-dependent change of variables defined in terms of the flow of the unperturbed system (\ref{eq:unpert}). The resulting integrator composes  $h/2$-flows of the unperturbed system (\ref{eq:unpert}) with RK steps applied to a transformed system, and we accordingly refer to them as flow-composed Runge-Kutta (FCRK)  methods.  Flow-composed Runge-Kutta exhibits (for a fixed integration interval) global errors of size $\mathcal{O}(\epsilon\, h^r)$ when applied to systems of the form (\ref{eq:pertODEeps}),  provided that the underlying RK scheme is of order $r$.  In addition, it is symplectic when both the perturbation and the unperturbed part are Hamiltonian and the underlying RK scheme is symplectic (in which case is necessarily implicit).

In the present work, we mainly focus on the application of flow-composed implicit Runge-Kutta (FCIRK)  methods with underlying 
implicit Runge-Kutta (IRK) schemes having favorable geometric properties. 
In practice, we choose as underlying IRK schemes the $s$-stage collocation methods with Gaussian quadrature nodes, which gives rise to time-symmetric symplectic schemes of the form 
 (\ref{eq:mStrang})  with global errors of size $\mathcal{O}(\epsilon\, h^{2s})$ when applied to perturbed systems of the form (\ref{eq:pertODEeps}).

  Compared to the application of IRK schemes to the original problem (\ref{eq:pertODEeps}), FCIRK method have several advantages: 
(i) while local truncation errors of size  $\mathcal{O}(h^{r+1})$ occur for a $r$th order IRK scheme,  a FCIRK scheme based on the same underlying IRK scheme has local truncation errors of size $\mathcal{O}(\epsilon\, h^{r+1})$; (ii) For implementations based on fixed point iterations, a faster convergence is expected for FCIRK schemes (the Lipschitz constant of the transformed system happens to be, roughly speaking, $\epsilon$ times smaller than that of the original perturbed system); (iii) FCIRK  methods admit mixed precision implementation (the IRK step for the transformed system implemented in a lower precision arithmetic, and the $h/2$-flows of the unperturbed system in higher precision) that allows us to efficiently reduce the effect of round-off errors.

In Section~\ref{s:FCRK} FCRK schemes are described in detail.  Structure-preserving properties of FCRK (mainly FCIRK) schemes are studied in Section~\ref{s:geo}. Implementation aspects of FCIRK methods are discussed in Section~\ref{s:impl}.
Section~\ref{s:ne} is devoted to present numerical experiments with a simplified N-body model of the solar system, where FCIRK schemes based on Gaussian nodes are compared to other methods. Some concluding remarks are given in Section~\ref{s:cr}.

\section{Flow-Composed Runge-Kutta  methods}
\label{s:FCRK}

Flow-Composed Runge-Kutta methods are based on an old idea, already exploited numerically by Lawson~\cite{Lawson1967} as early as  in 1967 for the case where $k(u)$ depends linearly on $u$. It consists on removing the relatively larger term $k(u)$  by applying a  time-dependent change of variables of the form 
\begin{equation}
\label{eq:globalchvar}
u = \varphi_{t-t_0}(U) 
\end{equation}
where $\varphi_{t}$ is the $t$-flow of (\ref{eq:unpert}). It is straightforward to check that in the new variables, the initial value problem (\ref{eq:pertODE})--(\ref{eq:iv}) reads
\begin{equation}
\label{eq:ODELawson}
\frac{d}{dt} U = \left(\varphi'_{t-t_0}(U)\right)^{-1} g\left(\varphi_{t-t_0}(U)\right), \quad U(t_0) = u_0,
\end{equation}
where  $\varphi'_{t}(U)$ represents the Jacobian matrix of $\varphi_{t}(U)$ with respect to $U$.   Lawson's generalized RK integration schemes 
can be interpreted as
numerically integrating  the initial value problem (\ref{eq:ODELawson}) by some Runge-Kutta scheme to  obtain approximations $U_j \approx U(t_j)$ of the solution $U(t)$ of (\ref{eq:ODELawson}) evaluated at $t=t_j$, $j=1,2,3,\ldots$, and then obtaining, back into the original variables,  approximations $u_j = \varphi_{t_j}(U_j) \approx u(t_j)$ of the solution   $u(t)$ of (\ref{eq:pertODE})--(\ref{eq:iv}).

The case where $k(u)= A u$, with $A$ a constant $D\times D$ matrix, have been extensively studied in the literature since the seminar work of~\cite{Lawson1967}. In that case,  $\varphi_{t}(u) = e^{t A} u$,  and the resulting methods belong to the broader family of exponential Runge-Kutta integrators~\cite{Hochbruck2010}, which are often referred to as Lawson RK methods~\cite{Celledoni2008,Cano2015},  or Integrating Factor RK methods~\cite{Bhatt2017}.

We have observed in some preliminary numerical experiments with perturbed Kepler problems, that the extension of Lawson's generalized RK methods to the case of nonlinear unperturbed system, only works well for relatively small integration intervals, but the performance of the method degrades (in particular, the local truncation errors become larger) for larger values of time $t$. 

Motivated by that, we propose computing  numerical approximations $u_j \approx u(t_j)$ of the solution   $u(t)$ of (\ref{eq:pertODE})--(\ref{eq:iv}) at times  $t=t_j$, $j=1,2,3,\ldots$ in a step-by-step manner as described next. The main idea is to make use  of a change of variables of the form
$u = \varphi_{t-\tau}(U)$, as in Lawson's method, but with a different value of the real parameter $\tau$ at each step.  

\begin{figure} [h!]
\centerline{\includegraphics [width=9cm, height=4cm] {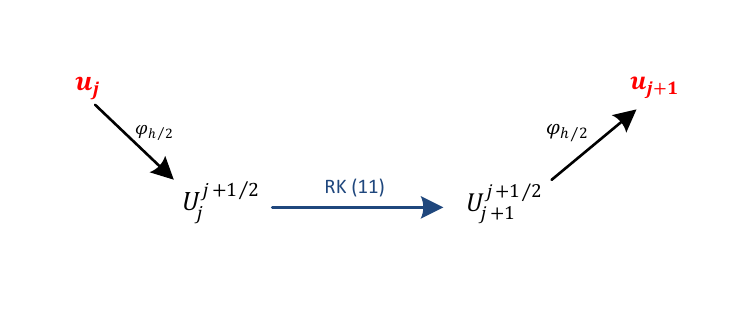}}
\caption{\small One step of a  FCRK method}
\label{fig:FCRKstep}
\end{figure}

We next give the definition of one step of a FCRK scheme. Assume that the approximation $u_{j} \approx u(t_{j})$ has already been computed, and that we want to compute an approximation  $u_{j+1} \approx u(t_{j+1})$, $t_{j+1}=t_j+h$, for the solution $u(t)$ of the initial value problem
\begin{equation}
\label{eq:pertODEj}
\begin{split}
&\frac{du}{dt} = k(u) + g(u), \\
&u(t_j) = u_j.
\end{split}
\end{equation}
The change of variables
\begin{equation}
\label{eq:localchvar}
u = \varphi_{t-\tau_j}(U^{j+1/2})
\end{equation}
%
where $\tau_j := \frac12(t_j + t_{j+1} )= t_j + \frac{h}{2}$,
transforms the problem (\ref{eq:pertODEj}) into
\begin{equation}
\label{eq:pertODEUj}
\begin{split}
&\frac{d}{dt} U^{j+1/2} = \left(\varphi'_{t-\tau_j}(U^{j+1/2})\right)^{-1} g\left(\varphi_{t-\tau_j}(U^{j+1/2})\right), \\
& U^{j+1/2}(t_j) = U^{j+1/2}_j,
\end{split}
\end{equation}
where $(t_j,U^{j+1/2}_j)$  and $(t_j,u_j)$ are related by the change of variables (\ref{eq:localchvar}), that is, $u_j= \varphi_{-h/2}(U^{j+1/2}_j)$. We thus propose applying one step of an RK scheme to approximately  compute $U^{j+1/2}(t_{j+1})$, and then approximating $u(t_{j+1})$ by undoing the change of variables.
To sum up, one step of a FCRK method is computed as follows: given $u_j \approx u(t_j)$ (see Figure~\ref{fig:FCRKstep}):
 \begin{enumerate}
 \item Compute 
 \begin{equation}
\label{eq:flow1}
U^{j+1/2}_j= \varphi_{h/2}(u_j).
\end{equation}
\item  Apply one step of length $h=t_{j+1}-t_{j}$ of the underlying RK method to obtain an approximation $U^{j+1/2}_{j+1}$
to the solution $U^{j+1/2}(t)$  of (\ref{eq:pertODEUj}) at $t=t_{j+1}$.
\item finally compute
\begin{equation}
\label{eq:flow2}
u_{j+1} = \varphi_{h/2}(U^{j+1/2}_{j+1}).
\end{equation}
 \end{enumerate}
%
%
%

%


\section{Geometric properties}
\label{s:geo}
We next discuss some geometric properties of FCRK methods when applied to autonomous perturbed systems of the form (\ref{eq:pertODE}).

\subsection{Affine changes of variables}

It is well known that applying an affine change of variables to a system of ODEs  and then numerically integrating the transformed system by means of an arbitrary  RK scheme, is equivalent to numerically integrating the original ODE with the same RK scheme and then applying the change of variables to translate the numerical result in terms of the new variables.

As a consequence of that, discretization by FCRK schemes of systems of the form (\ref{eq:pertODE}) and affine changes of variables also commute with each other.

It is worth mentioning that, in the particular case where $k(u)$ depends linearly on $u$, the numerical solution $u_j$ obtained by the application of a FCRK integrator coincides with the approximations obtained by the Lawson's generalized RK method based on the same underlying RK scheme. Indeed, if $k(u)=Au$, the new variables $U$ in Lawson's global change of variables (\ref{eq:globalchvar}) and the new variables $U^{j+1/2}$ in the local change of variables (\ref{eq:localchvar}) considered in FCRK integrators are related linearly as $U^{j+1/2} = e^{(j+\frac12) h A} U$. The fact that RK discretization and affine change of variables commute with each other implies that one step of a FCRK coincides with one step of Lawson's generalized RK method when $k(u)$ depends linearly on $u$.

\subsection{FCRK methods as one-step integrators}

We next show that each step of a FCRK scheme, summarized in Figure~\ref{fig:FCRKstep}, is of the form  (\ref{eq:mStrang}) when applied to autonomous systems of the form (\ref{eq:pertODE}).

Provided that $g(u,t)$ is independent of $t$, the approximation $U^{j+1/2}_{j+1}$
to the solution $U^{j+1/2}(t)$  of (\ref{eq:pertODEUj}) at $t=t_{j+1}$ can be equivalently seen as the approximation, 
 obtained by one step  of length $h$ of the underlying RK method,  to the value at $t=h$ of the solution $U(t)$
 of
 \begin{equation}
 \label{eq:pertODEU2}
\frac{d}{dt} U  = \left(\varphi'_{t-h/2}(U)\right)^{-1} g\left(\varphi_{t-h/2}(U) \right)
\end{equation}
with initial value $U(0) = U^{j+1/2}_j$.  The system (\ref{eq:pertODEU2}) is precisely obtained by applying the change of variables
\begin{equation}
\label{eq:chvar}
u=\varphi_{t-h/2}(U)
\end{equation}
to the original perturbed system (\ref{eq:pertODE}).

One step of a FCRK method can thus be defined as  (\ref{eq:mStrang}), where the map $\psi_h$ is defined  in terms 
of the real coefficients $a_{ij}$ ($1 \leq i, j \leq s$) and $b_{i}$, ($1\leq i \leq s$) of the underlying $s$-stage RK scheme as
\begin{equation}
\label{eq:psiFCIRK}
\psi_{h}(U) := U+ h \sum_{i=1}^s b_i \, f(c_i h, W_{i},h),
\end{equation}   
where 
\begin{equation}
\label{eq:f}
f(t,U,h):= \left(\varphi'_{t-h/2}(U)\right)^{-1} g\left(\varphi_{t-h/2}(U) \right),
\end{equation}
 $c_i = \sum_{j=1}^s a_{i j}$,  and  $W_{i}$ ($i=1,\ldots,s$) are defined (implicitly in the case of IRK methods) as functions of $(U,h) \in \mathbb{R}^{D+1}$ by
\begin{equation}
\label{eq:Ui}
U_{i} =U+ h \sum^s_{j=1}{a_{ij}\,f(c_j h,W_{j},h)}, \quad  i=1 ,\ldots, s.
\end{equation}

Flow-composed IRK schemes can be seen as generalizations of the second order integrator (\ref{eq:mStrang})--(\ref{eq:midpoint}): 
in the particular case where the underlying IRK scheme is chosen to be the implicit midpoint rule, the FCIRK method is reduced to the scheme (\ref{eq:mStrang})--(\ref{eq:midpoint}).

\subsection{Affine and quadratic invariants}

If the perturbed (\ref{eq:pertODE}) and the unperturbed (\ref{eq:unpert}) systems have a common first integral  $I(u)$, then it is straightforward to check that $I(u)$ is also an invariant of the problem (\ref{eq:pertODEU2}). 

If $I(u)$ is an affine invariant, since such kind of first integrals are conserved by arbitrary RK schemes, then $I(\psi_h(U)) \equiv I(U)$ for $\psi_h$ defined by (\ref{eq:psiFCIRK})--(\ref{eq:Ui}). Hence, $I(u)$ will be also preserved by the FCRK integrator (\ref{eq:mStrang}).

Similarly,  if $I(u)$ is a quadratic invariant and the underlying  $s$-stage IRK method is symplectic, that is, if 
\begin{equation} \label{eq:sympl_cond_1}
b_{i}a_{ij}+b_{j}a_{ji}-b_{i}b_{j}=0, \ \ 1 \leq i,j \leq s.
\end{equation}
then the value of  $I(u)$ will be exactly preserved by the FCIRK scheme.

\subsection{Application to Hamiltonian systems}

Let us assume that the system (\ref{eq:pertODE}) is of the form
\begin{equation}
\label{eq:unpertHam}
\frac{du}{dt} =  J^{-1} (\nabla K(u) + \nabla G(u)),
\end{equation}
where $K(u)$ and $G(u)$ are smooth real-valued functions, and $J$ is
the canonical skew-symmetric matrix 
\begin{equation}
\label{eq:J}
 J=\left(\begin{array}{cc}
   \ 0 & \ -I_d \\
     I_d & \ 0  \\
\end{array}\right)
\end{equation}
where $I_d$ is the identity matrix of size $d\times d$, $D= 2d$, and $u=(q,p)$, $q,p \in \mathbb{R}^d$. 

In that case, the $t$-flow $\varphi_{t}$ of the unperturbed system will be symplectic for each value of $t$. This implies that the system (\ref{eq:pertODEU2}) obtained by applying the change of variables (\ref{eq:chvar}) is a (non-autonomous) Hamiltonian ODE. The application of one step of length $h$ of a symplectic IRK scheme to (approximate the solution at $t=h$ of) a (non-autonomous or autonomous) Hamiltonian system with initial value at $t=0$ is a symplectic map as a function of the initial values. Hence, 
the map  (\ref{eq:mStrang}) (with $\psi_h$ given by (\ref{eq:psiFCIRK})--(\ref{eq:Ui})) corresponding to a step of a FCIRK method with symplectic underlying IRK scheme is also symplectic. 

It is well known that symplectic integrators are very appropriate for the long-term numerical integration of Hamiltonian systems~\cite{SanzSerna1994,Hairer2006}. The same good long-term behavior is also  guaranteed for conservative systems written in the more general form
\begin{equation}
\label{eq:pertODEHam2}
\frac{du}{dt} =  B\, (\nabla K(u) + \nabla G(u)),
\end{equation}
where $B$ is an arbitrary (constant) skew-symmetric matrix 
of size $D\times D$. It is always possible, up to a linear change of coordinates,  to write that system  in canonical coordinates so that
\begin{equation}
\label{eq:B}
B=\left(\begin{array}{ccc}
  \ 0 & \ I _d & 0\\
    -I_d & \ 0   & 0\\
     0  & 0 & 0
\end{array}\right)
\end{equation}
where $D \geq  2d$, and $u=(q,p,z)$, $q,p \in \mathbb{R}^d$, $z \in \mathbb{R}^{D-2d}$. The components of $z$ are casimirs of the Poisson structure given by the skew symmetric matrix $B$, that is, they are first integrals of the system (\ref{eq:pertODEHam2}) regardless of the actual definition of the Hamiltonian function $K(u)+G(u)$.

Clearly, numerically integrating by means of a FCRK scheme the system (\ref{eq:pertODEHam2}) in canonical coordinates (i.e., with $B$ given by (\ref{eq:B})) is equivalent to integrating with the same scheme the $2d$-dimensional Hamiltonian system obtained by replacing the components of $z$ by their constant values. Since FCRK discretization commute with linear changes of variables, it will also be equivalent to the numerical integration of 
(\ref{eq:pertODEHam2}) in non-canonical coordinates (i.e., with $B$ and arbitrary skew-symmetric matrix).

\subsection{Time-symmetry and reversibility}

A one-step method is said to be time-symmetric, if changing the sign of the time-step parameter $h$ in the definition of the one-step map of the method, has the effect of inverting that map. A IRK scheme with coefficients $a_{ij}$ ($1 \leq i, j \leq s$) and $b_{i}$, $c_{i}$ ($1\leq i \leq s$) is time-symmetric~\cite{Hairer2006} if (possibly after reordering the stages)
\begin{equation}
\label{eq:symmetry_cond}
\begin{split}
b_{s+1-i} &= b_{i}, \quad   c_{s+1-i} = 1 - c_{i}, \quad 1 \leqslant i \leqslant s,\\
b_{j} &= a_{s+1-i,s+1-j} + a_{i, j}, \quad 1 \leqslant i,j \leqslant s.
\end{split}
\end{equation}
It is straightforward to check that  the map $\psi_h$ given by (\ref{eq:psiFCIRK})--(\ref{eq:Ui})) satisfies
\begin{equation*}
\psi_{-h} = (\psi_{h})^{-1}.
\end{equation*}

This implies that a FCIRK scheme is time-symmetric, that is,
the map  (\ref{eq:mStrang}) (with $\psi_h$ given by (\ref{eq:psiFCIRK})--(\ref{eq:Ui})) corresponding to a step of a FCIRK method satisfies
\begin{equation*}
\widehat{\psi}_{-h} = (\widehat{\psi}_{h})^{-1},
\end{equation*}
if the underlying IRK scheme is time-symmetric.

Assume now that (\ref{eq:pertODE}) is $\rho$-reversible for a given invertible matrix $\rho \in \mathbb{R}^{D \times D}$, that is,
\begin{equation*}
k(\rho \, u) = -\rho \, k(u), \quad g(\rho\,  u) = -\rho\,  g(u).
\end{equation*}
Since FCIRK discretization commute with affine changes of variables, then the map  (\ref{eq:mStrang}) defining one-step of a FCIRK method satisfies that
\begin{equation*}
\widehat{\psi}_{h}(\rho\,  u) = \rho \, \widehat{\psi}_{-h}(u).
\end{equation*}
This implies that if the underlying IRK scheme is time-symmetric, then the FCIRK method is $\rho$-reversible in the sense that
\begin{equation*}
\rho^{-1} \widehat{\psi}_{h}(\rho \, u) = (\widehat{\psi}_{h})^{-1}(u).
\end{equation*}

%
%
%

\section{Implementation aspects of FCIRK integrators}
\label{s:impl}


Most scientific computing is done in double precision ($64$-bit) arithmetic but some problems need higher precision \cite{Bailey2012}. This class of problems can be solved by combining different precisions to achieve the required accuracy at minimal computational cost \cite{Dongarra2017}: most of the computation can be done in a lower precision (faster computation), while a smaller amount of computation can be done in higher precision (slower computation). Another way to increase the computing performance is based on the use of parallelism. FCIRK methods present several possibilities to parallelize the code.

\subsection{Mixed precision implementation}

Obtaining the numerical approximation  $u_{j} \approx u(t_{j})$, ($j=1,2,\ldots$) to the solution $u(t)$ of the initial value problem (\ref{eq:pertODE})--(\ref{eq:iv}) by means of an IRK scheme requires computing sums of the form 
\begin{equation}
\label{eq:sumu_j}
u_{j+1} = u_{j} + h\, f_j, \quad j=0,1,2,\ldots, 
\end{equation}
where each $f_j$ is a weighted sum of evaluations of the right-hand side of (\ref{eq:pertODE}).  
It is well known that round-off error accumulation in such sums can be greatly reduced with the use of Kahan's compensated summation algorithm~\cite{Kahan1965}. Compensated summation can be interpreted~\cite{Higham2002,Muller2009} from the point of view of backward error analysis of round-off errors as computing the exact sums of terms with relative errors comparable to the machine precision of the working floating point arithmetic. Computing the additions (\ref{eq:sumu_j}) with compensated summation is somehow equivalent to performing that sum in a floating point arithmetic with some additional extra digits of precision. 
The smaller the magnitude of the components of $h \, f_j$ relative to the components of $u_j$, the more additional digits of precision are gained with Kahan's compensated summation algorithm.
 
In  that case, the loss of precision due to round-off errors (in a given floating point arithmetic) in the application of an IRK scheme to a system of ODEs with right-hand side of smaller magnitude is reduced compared to the application of the same IRK scheme to an ODE with right-hand side of larger magnitude. 
  
This implies that, provided that  the perturbing terms in (\ref{eq:pertODE}) are of smaller magnitude than the terms of the unperturbed system (\ref{eq:unpert}), the errors due to round-off errors when computing $U^{j+1/2}_{j+1} = \psi_h(U^{j+1/2}_{j})$ (where the map $\psi_h$ is given by (\ref{eq:psiFCIRK})--(\ref{eq:Ui}))
) 
  in a FCIRK scheme will  be considerably smaller than in the application of one step of  the underlying IRK scheme to the original system (\ref{eq:pertODE}).  Of course, this advantage is lost if the applications of the $h/2$-flows in (\ref{eq:mStrang}) are implemented in the working precision arithmetic. This motivates us  to implement FCIRK  methods in mixed precision arithmetic: the IRK step for the transformed system is implemented in the working precision arithmetic, and the $h/2$-flows of the unperturbed system are implemented in a higher precision one.  This technique allows us to efficiently reduce the effect of round-off errors.

\subsection{Efficient evaluation of the right-hand side of  the transformed ODE system}

Each evaluation of the right-hand side of that transformed system (\ref{eq:pertODEU2}) requires: 
\begin{itemize}
\item  computing $R=g(\varphi_{t-h/2}(U))$, 
\item  followed by the computation of $(\varphi'_{t-h/2}(U) )^{-1} R$.
\end{itemize}
 The later  can be efficiently computed provided that the unperturbed system (\ref{eq:unpert}) is Hamiltonian.  
Indeed, let us assume that the unperturbed system (\ref{eq:unpert}) is of the form
\begin{equation}
\label{eq:unpertHam}
\frac{du}{dt} =  J^{-1} \nabla K(u),
\end{equation}
where $K(u)$ is a smooth real-valued function, and $J$ is 
the canonical skew-symmetric matrix (\ref{eq:J}). 
By symplecticity of the $t$-flow $\varphi_{t}$ of (\ref{eq:unpertHam}),  
\begin{equation*}
(\varphi'_{t}(u))^T J \varphi'_{t}(u)= J, \quad 
\forall(t,u) \in \mathbb{R}^{D+1},
\end{equation*}
and hence, for arbitrary $R \in \mathbb{R}^D$
\begin{equation*}
(\varphi'_{t}(u) )^{-1} R = J^{-1}(\varphi'_{t}(u))^{T} J R.
\end{equation*}
Hence, the right hand side of (\ref{eq:pertODEU2}) can be evaluated as follows: 

\begin{algorithm}[H]
 \BlankLine
 $u=\varphi_{t-h/2}(U)$;
 \BlankLine
 $R=g(u)$; $\widehat R=J R$;
 \BlankLine
 $\widehat F = \varphi'_{t-h/2}(U)^{T} \widehat R$;
 \BlankLine 
 $F = J^{-1} \widehat F$ ;
 \BlankLine
\end{algorithm}
It is important to observe that  $\varphi'_{t-h/2}(U)^T \widehat R$ is the gradient with respect to $U$ of the scalar function $\varphi_{t-h/2}(U)^T \widehat R$. 
Hence, the vector  $\widehat F =\varphi'_{t-h/2}(U)^T \widehat R$ in the algorithm above  can be  computed very efficiently  by using the technique of reverse mode automatic differentiation~\cite{Linnainmaa1976,Griewank2008,Griewank2012}.  In that technique, the intermediate results needed in the computation of $\varphi_{t-h/2}(U)$ are reused for an efficient computation of the gradient of  $\varphi_{t-h/2}(U)^T \widehat R$.

\subsection{Parallelization} 

When computing one step of the underlying  $s$-stage IRK method applied to (\ref{eq:pertODEU2}), the right-hand side of that transformed system needs to be evaluated  $s$ times at each fixed point iteration.
Clearly, these $s$  evaluations can be performed in parallel if $s$ processors or cores are available. 

In addition, if the unperturbed system consists of a collection of uncoupled subsystems of ODEs (as in the solar system model considered in Section~\ref{s:ne}), then the applications of the $h/2$-flows in Figure~\ref{fig:FCRKstep} consists of several uncoupled flows that can also be computed in parallel.

\subsection{Some additional implementation details}

If output is not required at each step,  some computational saving is obtained by composing two consecutive $h/2$-flows of the unperturbed system: if $u_{j}$ does not need to be computed, one can reduce the sequence of computations 
\begin{equation*}
u_{j} = \varphi_{h/2}(U^{j-1/2}_{j}), \quad U^{j+1/2}_{j} = \varphi_{h/2}(u_j), 
\end{equation*}
into just one flow evaluation $U^{j+1/2}_{j}=\varphi_{h}(U^{j-1/2}_{j})$.

In our implementation used to make the numerical experiments presented in next section, we actually organize the computations as follows provided that output is required every $m$ steps:

 \begin{algorithm}[H]
  \BlankLine
  $U^{1/2}_{0} =\varphi_{h/2}(u_{0})$\;
   \BlankLine
   \For{$j=1,2,3,\ldots$}
   {
    \BlankLine
     $U^{j-1/2}_{j}=\psi_{h}(U^{j-1/2}_{j-1})$\;
   $U^{j+1/2}_{j}=\varphi_{h}(U^{j-1/2}_{j})$\;
    \If{$\mathrm{mod}(j,m)=0$}
    {$u_j = \varphi_{-h/2}(U^{j+1/2}_{j})$\;
    }
    \BlankLine
   }
 \end{algorithm}

\subsection{Application to non-autonomous problems}

It is worth mentioning that FCRK methods can be applied to non-autonomous perturbed systems where the unperturbed and the perturbing term  depend explicitly on $t$.  The resulting scheme for non-autonomous problems can be derived from the corresponding integration method for autonomous problems in the standard way: it is enough to consider the application of the method  to the formally autonomous problem  obtained by adding the time variable $t$ to the set of state variables and considering the additional differential equation obtained by setting the derivative of $t$ as 1.

We next explicitly formulate the resulting method 
in the particular case of non-autonomous systems of the form
\begin{align}
\label{eq:pertODEt}
&\frac{du}{dt} = k(u) + g(t,u), \\ 
\label{eq:ivt}
&u(t_0) = u_0 \in \mathbb{R}^D.
\end{align}
where only the perturbing term depends explicitly on time.

Numerical approximations $u_j \approx u(t_j)$ of the solution  $u(t)$ of (\ref{eq:pertODEt})--(\ref{eq:ivt}) at times  $t=t_j$, $j=1,2,3,\ldots$ are computed by a FCRK method as follows: given $u_j \approx u(t_j)$,
 \begin{enumerate}
 \item compute 
 \begin{equation}
\label{eq:flow1}
U^{j+1/2}_j= \varphi_{h/2}(u_j).
\end{equation}
\item  compute 
\begin{equation}
\label{eq:psiFCIRKt}
U^{j+1/2}_{j+1} = U^{j+1/2}_{j}+ h \sum_{i=1}^s b_i \, f_j(c_i h, W_{j,i},h),
\end{equation}   
where 
\begin{equation}
\label{eq:fj}
f_j(t,U,h):= \left(\varphi'_{t-t_j-h/2}(U)\right)^{-1} g\left(t,\varphi_{t-t_j-h/2}(U) \right),
\end{equation}
 $c_i = \sum_{j=1}^s a_{i j}$,  and  $W_{j,i} \in \mathbb{R}^D$ ($i=1,\ldots,s$) are determined (implicitly in the case of IRK methods)  by
\begin{equation}
\label{eq:Wji}
W_{j,i} =U^{j+1/2}_{j}+ h \sum^s_{\ell=1}{a_{i\ell}\,f_j(c_{\ell} h,W_{j,\ell},h)}, \quad  i=1 ,\ldots, s.
\end{equation}
\item finally compute
\begin{equation}
\label{eq:flow2}
u_{j+1} = \varphi_{h/2}(U^{j+1/2}_{j+1}).
\end{equation}
 \end{enumerate}

\section{Numerical experiments}
\label{s:ne}

We next report some numerical experiments with the $10$-body model of the solar system written in canonical heliocentric coordinates. 
We have compared FCIRK methods with high order composition and splitting methods. Efficiency diagrams of the comparison are given for different precision of the FCIRK integrators: one fully implemented in double precision ($64$-bit), second one with a "double - long double" mixed precision implementation and third one with a "long double - quadruple" mixed precision implementation. 
Finally, we have presented an statistical analysis of the propagation of round-off errors.

\subsection{Test problem: $10$-body model of the solar system}

As test problem, we consider a Newtonian 10-body model of the solar system, with point masses corresponding to the sun, the eight planets and Pluto (with the barycenter of the Earth-Moon system instead of the position of the center of the Earth). In all of the numerical experiments presented here we consider an integration interval of $10^6$ days.

The equations of motion are
\begin{equation*}
\frac{d^2}{dt^2} q_i= G \sum_{j=0,j \neq i}^{N} \frac{m_j}{\|q_j-q_i\|^3} (q_j-q_i) , \ \  i=0,1,\dots, N,
\end{equation*}
where $N=9$,   $q_i\in \mathbb{R}^3$, $m_i \in \mathbb{R}, \ \ i=0,\dots,N$ are the positions and mass of each body respectively, $G$ is the gravitational constant, and $\|\cdot\|$ represents the Euclidean norm. We have considered $m_0$ and $q_0$ to be the mass and position of the sun.

This is a Hamiltonian system with Hamiltonian function
\begin{equation*}
H(q,p)=\frac{1}{2}\ \sum^N_{i=0}{\ \frac{{\|p_i\|}^2}{m_i}}-G\ \sum^N_{0\le i<j\le N}{\frac{m_im_j}{\|q_i-q_j\|}},
\end{equation*}
where $p_i = m_i \dot{q}_i$, $i=0,1,\ldots,N$. We consider initial values taken from~\cite{Folkner2014}, they are shifted so that the barycenter is at the origin of coordinates and such that the system has vanishing linear momentum,
\begin{equation*}
\sum_{i=0}^N m_i \, q_i(0)=0, \quad \sum_{i=0}^N p_i(0)=0.
\end{equation*}
The equations of motion thus gives the time evolution of the barycentric coordinates.

In order to integrate the system by means of FCIRK methods, we  rewrite the equations of motion as a perturbation of 9 uncoupled Keplerian systems. This can be accomplished by rewriting the original system in Jacobi coordinates, or in canonical heliocentric coordinates. Here, we choose the second option. Canonical heliocentric coordinates $Q_i,P_i \in \mathbb{R}^3, \ i=0,\dots,N$ are defined as follows:
\begin{align*}
Q_0 &=\frac{1}{M}\sum_{i=0}^{N} m_i \, q_i, \quad P_0 =\sum_{i=0}^{N}p_i, \\
Q_i &=q_i-q_0, \quad  P_i=p_i - \frac{m_i}{M} P_0, \quad  i=1,\dots,N,
\end{align*}
where $M=\sum_{i=0}^{N} m_i$. 

The Hamiltonian function in the canonical heliocentric coordinates reads
\begin{equation*}
H(Q,P):=\frac{1}{2m_0} \|P_0\|^2 + H_K(Q,P) + H_I(Q,P),
\end{equation*}
where
\begin{align}
\label{eq:HK}
 &H_K(Q,P) := \sum_{i=1}^{N}\bigg(\frac{\|P_i\|^2}{2 \mu_i} - \mu_i \frac{k_i}{\|Q_i\|}\bigg),\\
 \label{eq:HI}
 &H_I(Q,P) := \frac{1}{m_0} \left(\sum\limits_{1 \leqslant i<j}^{N} P_i^T\ P_j \right)  -\sum\limits_{1 \leqslant i<j}^{N} \frac{G m_i m_j}{\|Q_i-Q_j\|},
\end{align}
with
\begin{equation*}
\frac{1}{\mu_i} = \frac{1}{m_i} + \frac{1}{m_0}, \quad 
k_i = G \, (m_0 +m_i), \quad 
i=1,\ldots,N. 
\end{equation*}
The equations of motion for $Q_0, P_0$ may be ignored, as $\frac{d}{dt} P_0=0$, which implies that if $q_i,p_i$ are barycentric coordinates, then the center of mass $Q_0$ remains at the origin of coordinates with zero linear momentum $P_0$.

Instead of writing the equations of motion for $Q_i,P_i$, $i=1,\ldots,N$, we find it more convenient to consider the system of ODEs corresponding to the variables $Q_i$ and $V_i:=P_i/\mu_i$. We thus have the system
\begin{equation}
\label{eq:EDOQV}
\frac{d}{dt} 
\left(
\begin{matrix}
Q_i \\
\\
V_i
\end{matrix}
\right) =
\left(
\begin{matrix}
V_i \\
\\
\displaystyle
- \frac{k_i}{\|Q_i\|^3 }\ Q_i
\end{matrix}
\right) +
\sum_{j\ne i,\ j=1}^{N} 
g_{i,j}(Q,V), \quad i=1,\ldots,N,
\end{equation}
where
\begin{equation*}
g_{i,j}(Q,V) =
\left(
\begin{matrix}
\displaystyle
 \frac{V_j \ m_j}{(m_0+m_j)} \\
 \\
 \displaystyle
-\frac{k_i}{m_0} \ \frac{m_j}{\|Q_i-Q_j\|^3} (Q_i-Q_j)  
\end{matrix}
\right).
\end{equation*}
This is a system of the form (\ref{eq:pertODEHam2}), with 
\begin{equation*}
B=\left(\begin{array}{cc}
   \ 0 & \mathcal{M}^{-1} \\
     -\mathcal{M}^{-1} & \ 0  \\
\end{array}\right),
\end{equation*}
where $\mathcal{M}$ is a diagonal matrix with entries
\begin{equation*}
(\mu_1,\mu_1,\mu_1,\mu_2,\mu_2,\mu_2,\ldots,\mu_N,\mu_N,\mu_N),
\end{equation*}
and
\begin{align*}
K(Q,V) &= \sum_{i=1}^{N} \mu_i \, \bigg(\frac{\|V_i\|^2}{2} -  \frac{k_i}{\|Q_i\|}\bigg), \\
G(Q,V) &= \frac{1}{m_0} \left(\sum_{1 \leqslant i<j}^{N} \mu_i \ \mu_j V_i^T\ V_j \right)  -\sum_{1 \leqslant i<j}^{N} \frac{G m_i m_j}{\|Q_i-Q_j\|}.
\end{align*}

\subsection{Considered numerical integrators}

We will compare FCIRK methods based on IRK collocation schemes with Gauss-Legendre nodes to high order explicit symplectic integrators.

The generalized leapfrog method is commonly known in the context dynamical astronomy as  Wisdom and Holman's symplectic map~\cite{Wisdom1991}. With Poincare's canonical heliocentric coordinates, the $t$-flow of the perturbing Hamiltonian cannot be explicitly computed. Of course, the integration scheme (\ref{eq:mStrang})--(\ref{eq:midpoint}) could be applied,  but it is not an explicit integrator.  An explicit alternative is usually preferred in that case. Indeed, the perturbation of the N-body problem written in canonical heliocentric coordinates  is a sum of two terms (\ref{eq:HI}), one depending only on the generalized momenta, and the other one depending only on the generalized positions. In that case,  the implicit midpoint rule (\ref{eq:midpoint}) can be effectively replaced in (\ref{eq:mStrang}) by the application of Strang splitting to the perturbation itself split in the position-depending and momentum-depending  parts~\cite{Touma1993,Farres2013}. This second order integrator is however not very well suited for high accuracy computations.
 
 Explicit symplectic methods of higher order of accuracy based on splitting for perturbed problems (\ref{eq:pertODEeps}) are constructed in~\cite{Blanes2013}, and tested for simulations of simple N-body models of  the solar system in~\cite{Farres2013}.  Methods designed for its application to perturbed systems where the flow of the perturbation is replaced by a second order discrete approximation are constructed and analyzed as well. Among them, the method \emph{ABAH1064}, with errors of size  $\mathcal{O}(\epsilon\, h^{10} + \epsilon^2\ h^6 + \epsilon^3\ h^4)$ when applied to a perturbed system of the form (\ref{eq:pertODEeps}), showed excellent performance when applied to solar system simulations in canonical heliocentric coordinates. 
 
 We will mainly compare the performance of our FCIRK method to the explicit integrator \emph{ABAH1064}. We will also consider, as a reference, an explicit symplectic Runge-Kutta-Nystr{\o}m method obtained as a composition method based on the standard leapfrog method. We consider the $10$th order composition scheme with $35$ stages of Sofroniu and Spaletta~\cite{Sofroniou2005} (see also \cite{Hairer2006}) which we refer as \emph{CO1035}.

Compared to explicit symplectic methods, FCIRK schemes have the disadvantage of having to solve a set of implicit equations at each step.  We accomplish this by making use of the implementations of symplectic  IRK methods  based on fixed point iteration and on Newton-like iterations presented in~ \cite{antonana2017}  and \cite{antonana2017a} respectively. We have observed that,  for solar system simulations, the implementation based on fixed point iterations is more efficient than the implementation based on Newton iterations, so that we restrict ourselves to the former.


\subsection{Efficiency comparisons}

All the computations presented here are made on Intel Xeon processor E5-2640 ($8$ cores) with Ubuntu $14.04$ Linux operating system and we have used \emph{gcc} compiler.

Three kinds of efficiency diagrams are employed in the present subsection. In the first one, the number of evaluations of the perturbation terms versus the maximum error in energy  in double logarithmic scale is displayed. In the second one, the cpu time for sequential executions is considered instead of the number of evaluations. In the third kind of diagram, wall time is displayed with execution of the implicit schemes (FCIRK) made with some parallelization. More precisely, the $s$ evaluations of the right-hand side required at each iteration are performed with OpenMP using $4$ threads. No parallelization is exploited for the evaluation of the flows of the unperturbed system.

\subsubsection{FCIRK and IRK comparison}

In the first experiment we compare FCIRK methods with the application of its underlying IRK scheme to the original system (\ref{eq:pertODE}). We have compared the performance of the $12$th order $6$-stage IRK collocation method  based on the IRK Gauss-Legendre nodes with the corresponding FCIRK integrator.

For the IRK scheme, we consider the partitioned fixed point iteration considered in \cite{Hairer2006}, and carefully implemented in~\cite{antonana2017} in double precision arithmetic. Among the several initializations of the fixed point iterations implemented in~\cite{antonana2017},  the best choice for this problem seems to be the initialization based on the  interpolation of previous stages. 

\begin{figure}[h!]
\centering
\begin{tabular}{c c}
\subfloat[]
{\includegraphics[width=.45\textwidth]{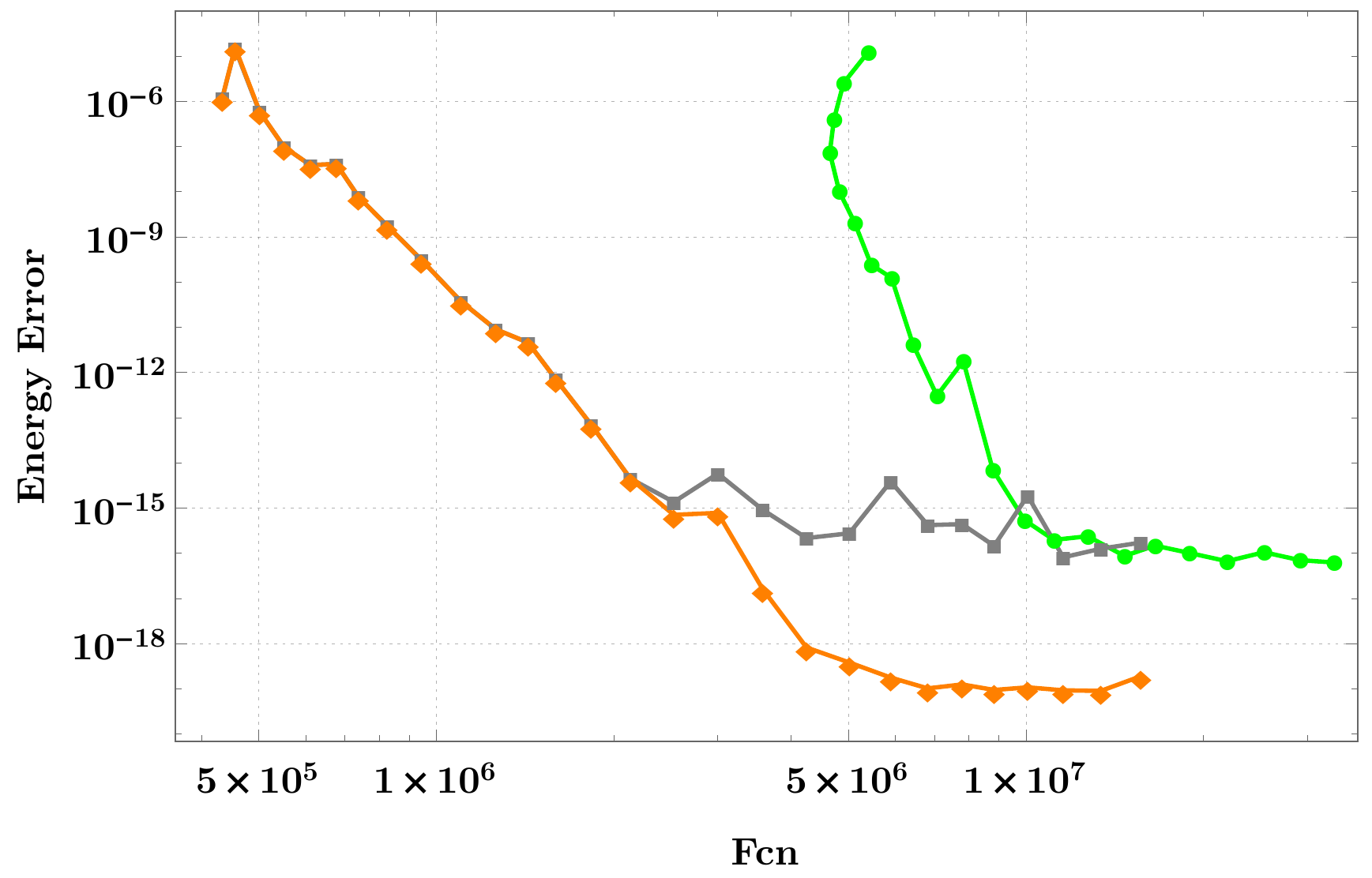}}
&
\subfloat[]
{\includegraphics[width=.45\textwidth]{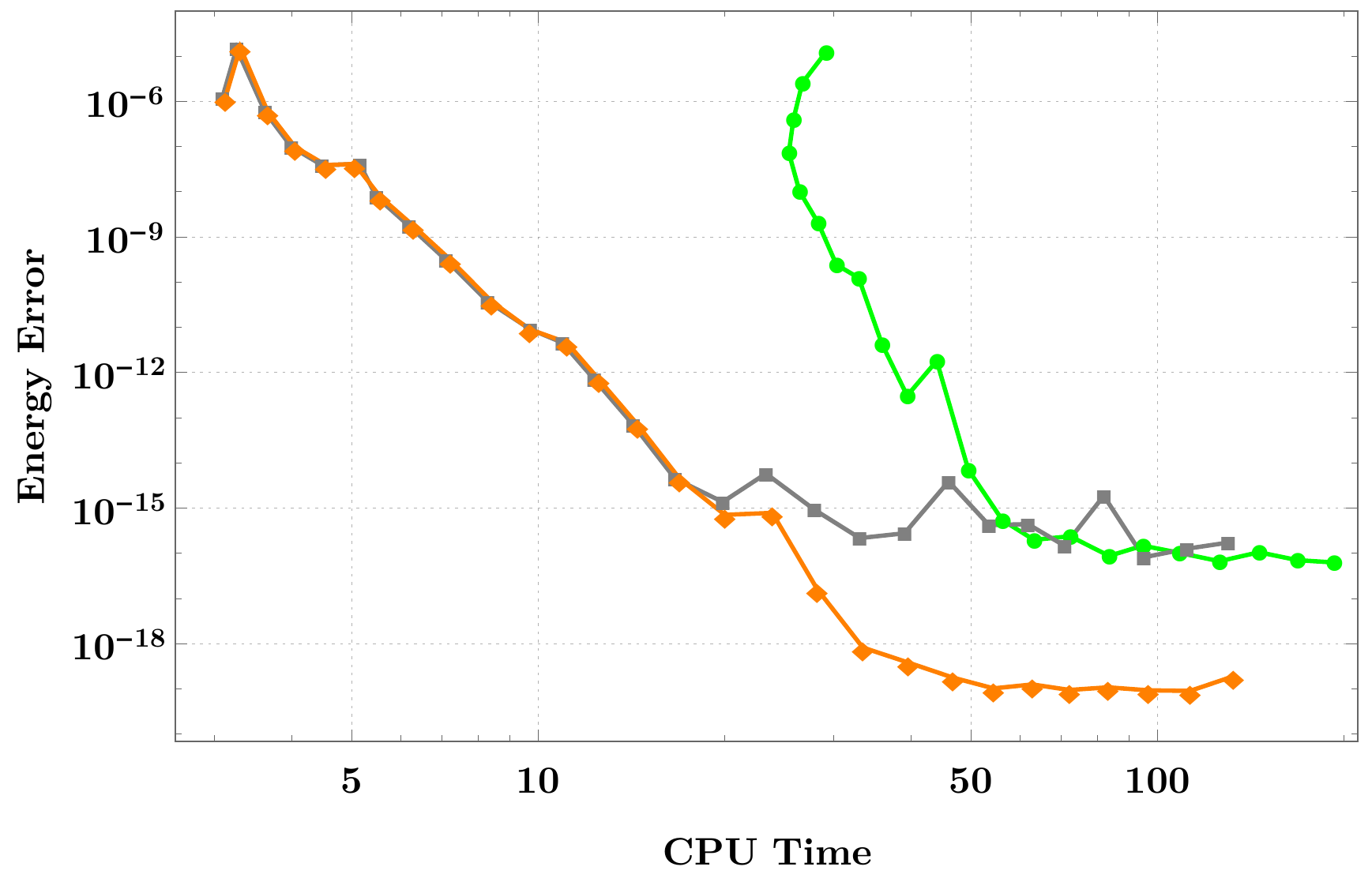}}
\end{tabular}
\caption{\small
 Efficiency diagrams (in double logarithmic scale) of methods based on 6-stage collocation methods of order $12$: IRK double precision arithmetic (green), FCIRK in double precision arithmetic (gray), and FCIRK in mixed precision (double - long double) arithmetic (orange): (a) Number of function evaluations versus  error in energy, (b) cpu time versus error in energy}
\label{fig:esp1}
\end{figure}

In Figure~\ref{fig:esp1} the efficiency diagrams for the two integrators are compared. Two different versions of the FCIRK integrators are considered. One fully implemented in double precision ($64$-bit) arithmetic, and the other one with a mixed-precision implementation: long double ($80$-bit) precision for the evaluations of flows in (\ref{eq:mStrang}), and double precision arithmetic for the application of the underlying IRK scheme to the transformed problem. In the first plot, where the number of function evaluations versus error in energy is shown, the advantage of the FCIRK method over the IRK scheme is more pronounced than in the second plot, where cpu time versus error in energy is displayed. The difference among the two efficiency measures accounts for the relative overhead of computing the flows of the unperturbed system (and related computations) required in the FCIRK problem. 
The superiority of the FCIRK integrator is clear in any case.

%

\subsubsection{FCIRK $s=6,8,16$ methods}

\begin{figure} [h!]
\centerline{\includegraphics [width=8cm, height=5cm] {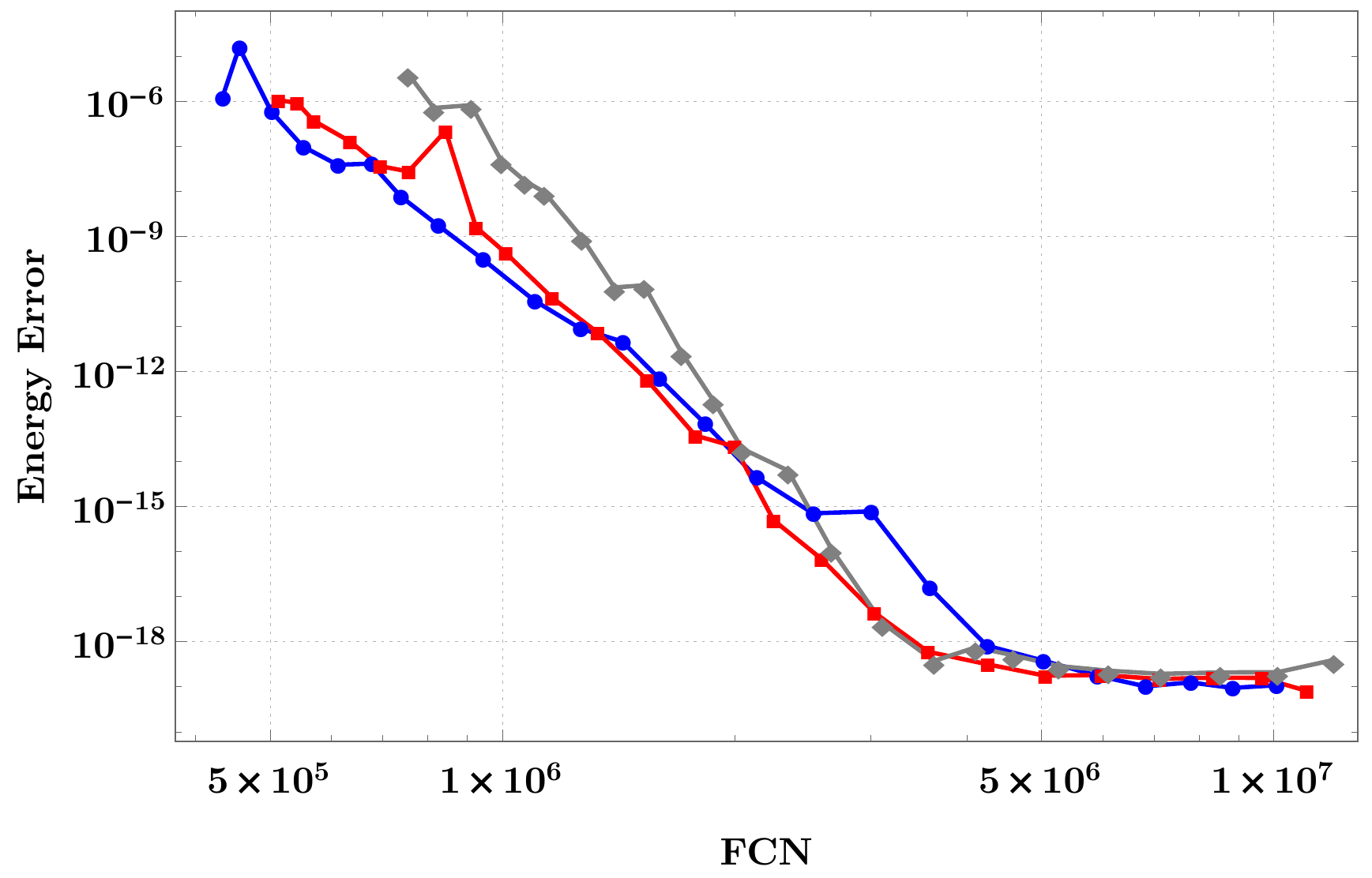}}
\caption{\small Efficiency diagrams (in double logarithmic scale) of the number of evaluations of the perturbations terms versus the maximum error in energy. We compare three FCIRK methods implemented with mixed precision (double - long double): methods with $s=6$ stages (blue), $s=8$ stages (red) and $s=16$ stages (gray)}
\label{fig:esp2a}
\end{figure}

\begin{figure}[h!]
\centering
\begin{tabular}{c c}
\subfloat[Sequential executions]
{\includegraphics[width=.45\textwidth]{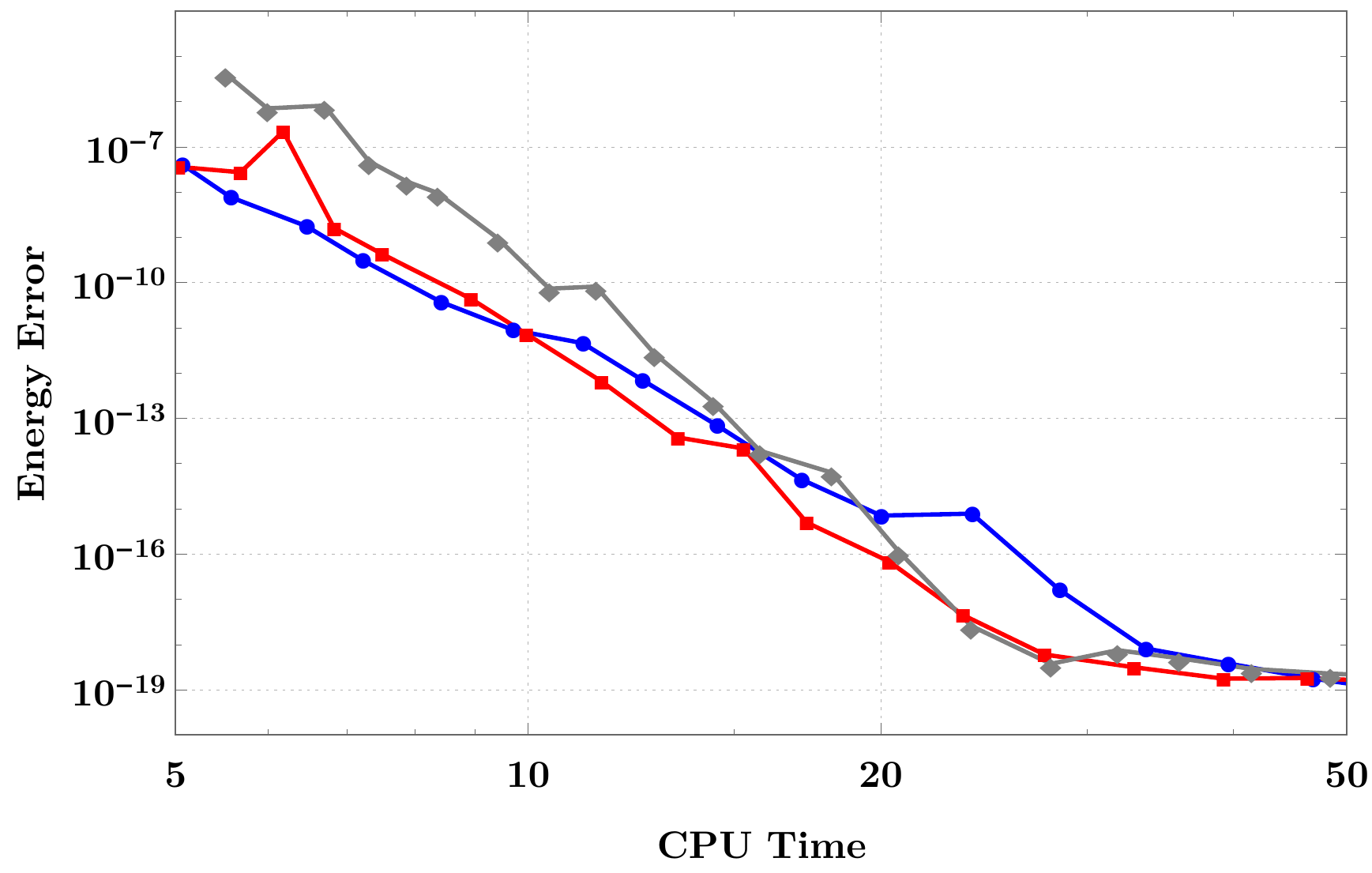}}
&
\subfloat[Parallel executions]
{\includegraphics[width=.45\textwidth]{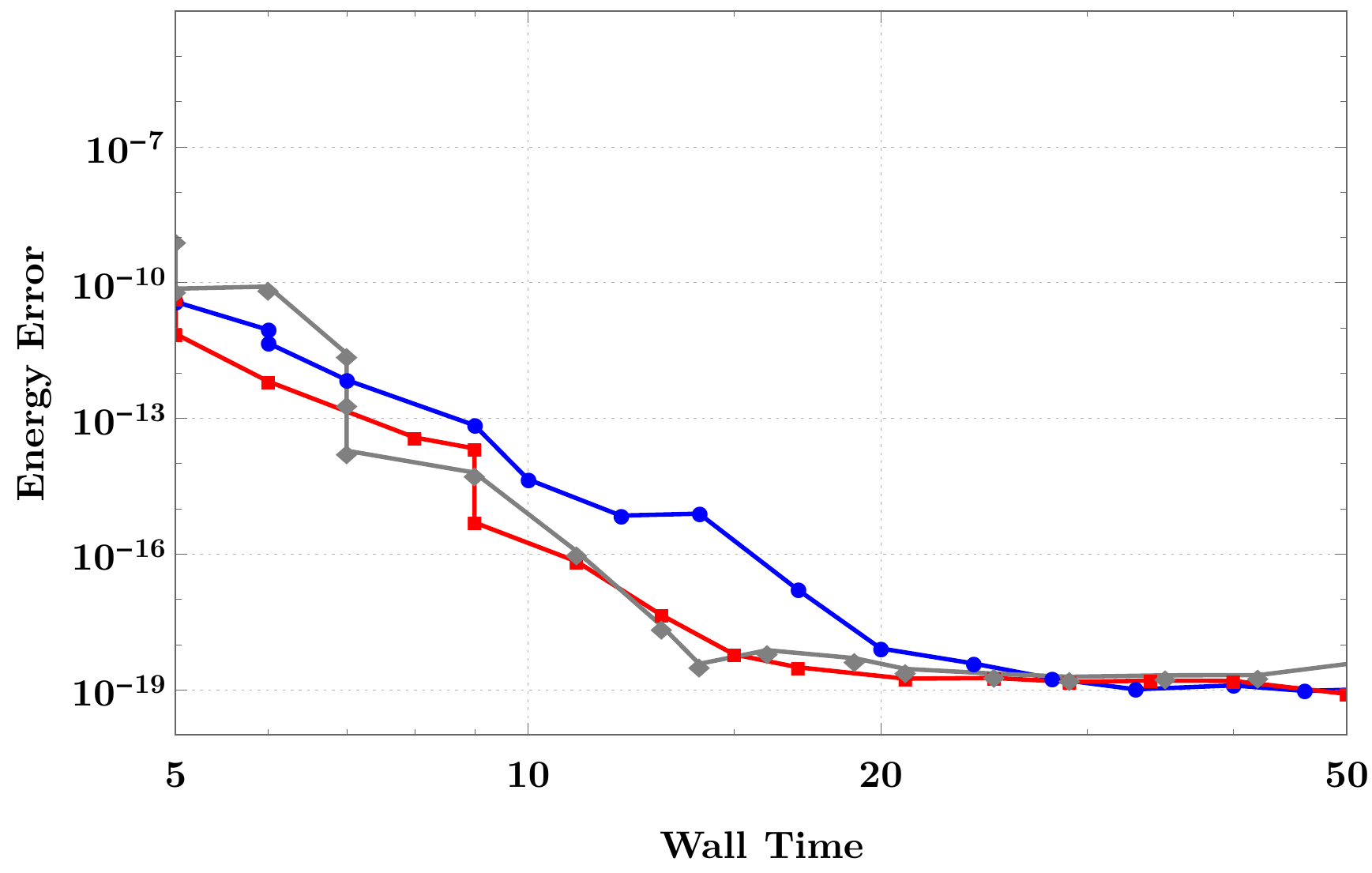}}
\end{tabular}
\caption{\small
Efficiency diagrams (in double logarithmic scale) of (a) sequential executions cpu time  and (b) parallel executions wall time using threads=$4$ versus the maximum error in energy. We compare three FCIRK methods implemented with mixed precision (double - long double): methods with $s=6$ stages (blue), $s=8$ stages (red) and $s=16$ stages (gray)}
\label{fig:esp2b}
\end{figure}

In Figures~\ref{fig:esp2a} and \ref{fig:esp2b} the efficiency of FCIRK methods based on IRK Gauss collocation methods of $s=6,8,16$ stages implemented with mixed precision  (double - long double) arithmetic is compared. For higher accuracy computations the methods with $s=8$ and $s=16$ are clearly superior to the method with $s=6$.

\subsubsection{FCIRK, \emph{ABAH1064} and \emph{CO1035} comparison}

\begin{figure} [h!]
\centerline{\includegraphics [width=8cm, height=5cm] {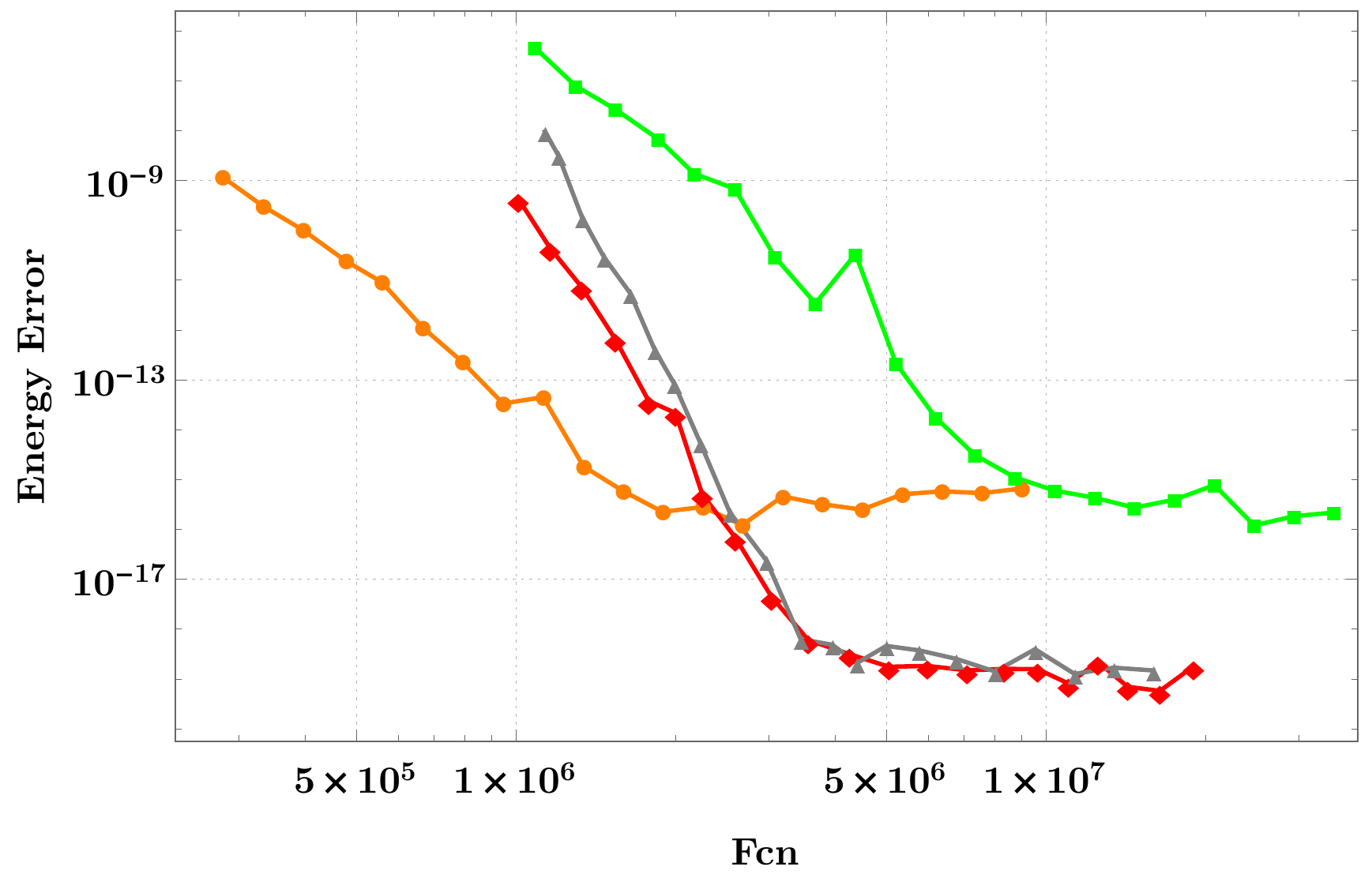}}
\caption {\small Efficiency diagrams (in double logarithmic scale) of the number of evaluations of the perturbations terms versus the maximum error in energy. We compare FCIRK methods with $s=8$ stages (red) and $s=16$ stages (gray) implemented with mixed precision (double - long double) and explicit integrators \emph{ABAH1064} (orange) and \emph{CO1035} (green) implemented in double precision arithmetic}
\label{fig:esp3a}
\end{figure}

\begin{figure}[h!]
\centering
\begin{tabular}{c c}
\subfloat[Sequential executions]
{\includegraphics[width=.45\textwidth]{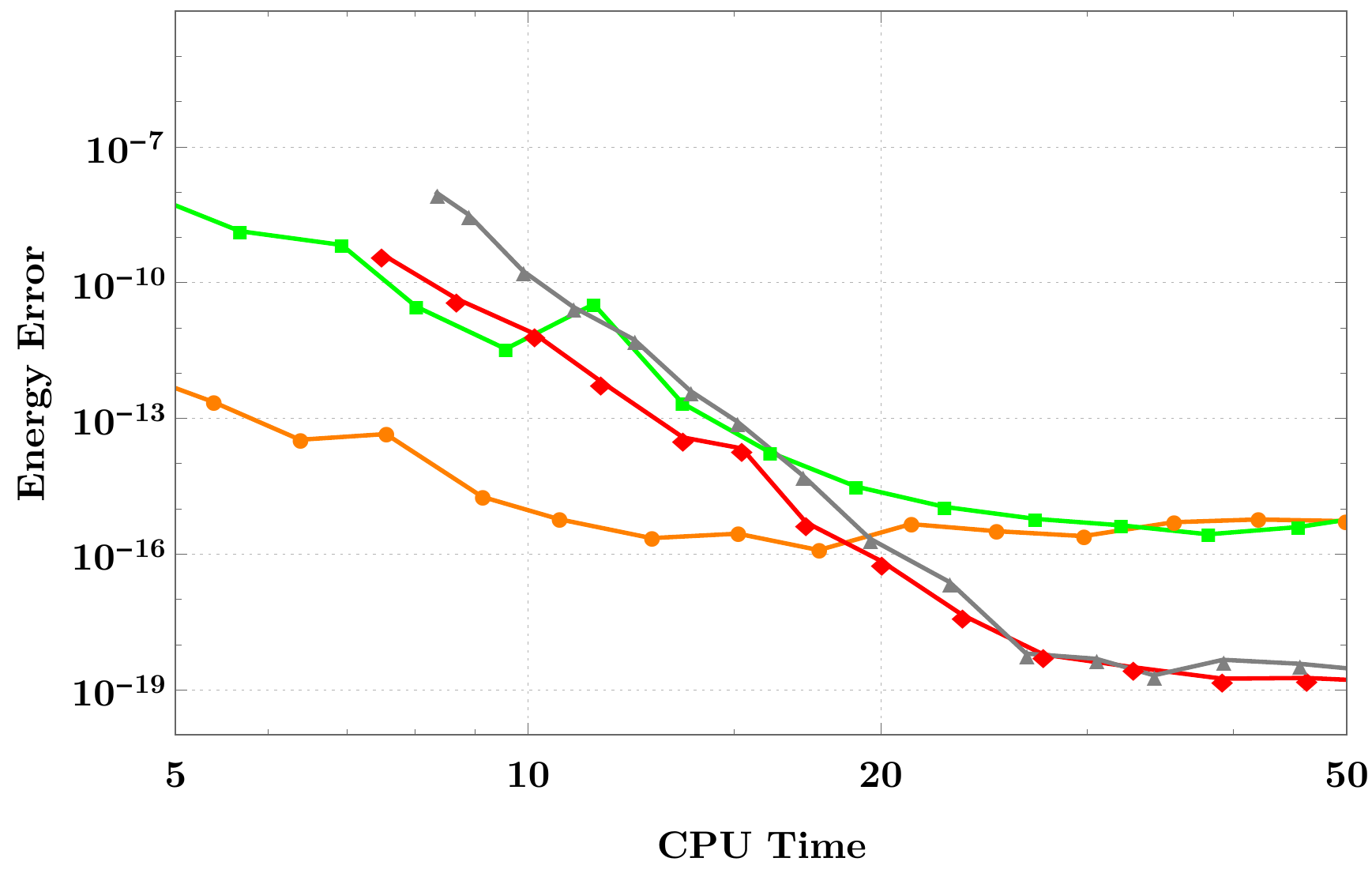}}
&
\subfloat[Parallel executions]
{\includegraphics[width=.45\textwidth]{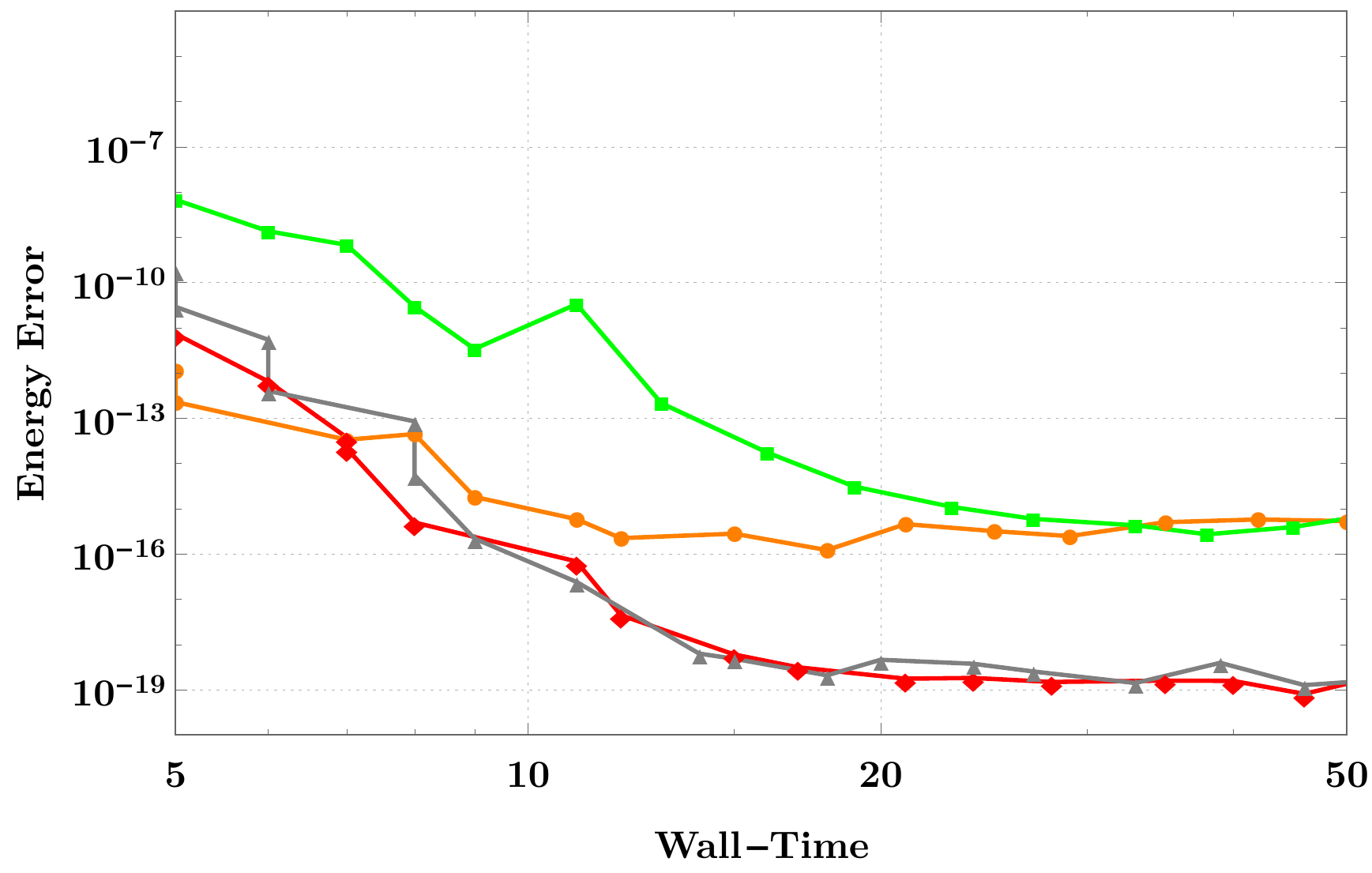}}
\end{tabular}
\caption {\small
Efficiency diagrams (in double logarithmic scale) of (a) sequential executions cpu time  and (b) parallel executions wall time using threads=$4$ versus the maximum error in energy. We compare FCIRK methods with $s=8$ stages (red) and $s=16$ stages (gray) implemented with mixed precision (double - long double) and explicit integrators \emph{ABAH1064} (orange) and \emph{CO1035} (green) implemented in double precision arithmetic}
\label{fig:esp3b}
\end{figure}

In Figures~\ref{fig:esp3a} and \ref{fig:esp3b} the efficiency of FCIRK methods based on IRK Gauss collocation methods of $s=8,16$ stages implemented with mixed precision  (double - long double) arithmetic is compared to the explicit integrators   \emph{ABAH1064} and \emph{CO1035} implemented in double precision arithmetic. For sequential executions,  \emph{ABAH1064} is more efficient than the FCIRK methods, up to the highest precision allowed by its implementation in double precision arithmetic. The simple parallelization technique implemented for the FCIRK methods makes possible to achieve execution times similar to those of  \emph{ABAH1064}.

When higher accuracy simulations are required, the effects of round-off errors should be reduced. Performing all the computations in quadruple precision arithmetic is very expensive in computing time.  In solar system simulations with explicit symplectic integrators, the influence of round-off error  is sometimes minimized by using $80$-bit precision (long double) arithmetic for all the computations~\cite{Laskar2011,Laskar2015}. Hence, in this case we have computed explicit integrators \emph{ABAH1064} and \emph{CO1035} using long double ($80$-bit) precision arithmetic. 

\begin{figure} [h!]
\centerline{\includegraphics [width=8cm, height=5cm] {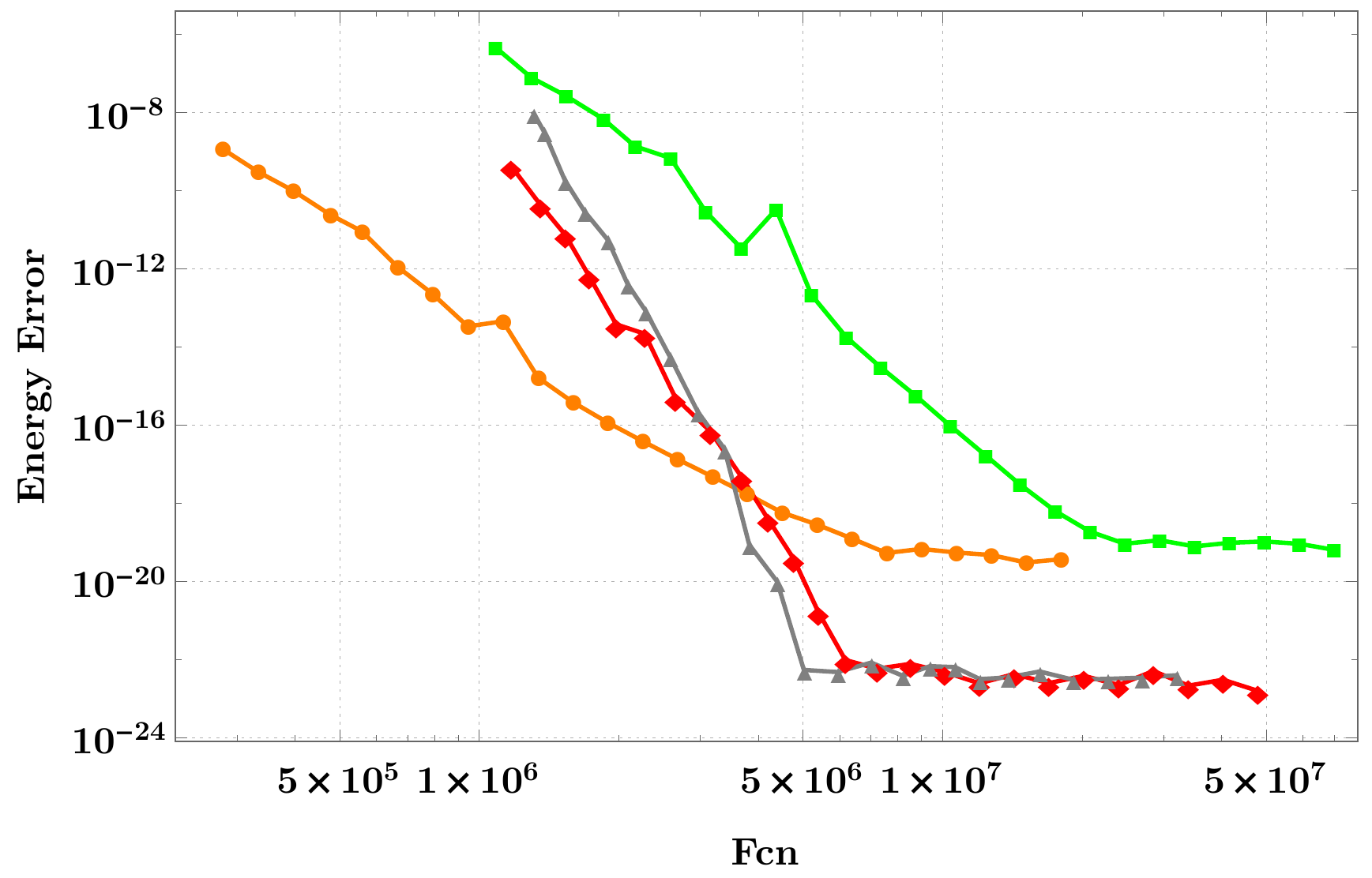}}
\caption{\small Efficiency diagrams (in double logarithmic scale) of the number of evaluations of the perturbations terms versus the maximum error in energy. We compare FCIRK methods with $s=8$ stages (red) and $s=16$ stages (gray) implemented with mixed precision (long double - quadruple) and explicit integrators \emph{ABAH1064} (orange) and \emph{CO1035} (green) implemented in long double ($80$-bit ) precision arithmetic}
\label{fig:esp4a}
\end{figure}

\begin{figure}[h!]
\centering
\begin{tabular}{c c}
\subfloat[Sequential executions]
{\includegraphics[width=.45\textwidth]{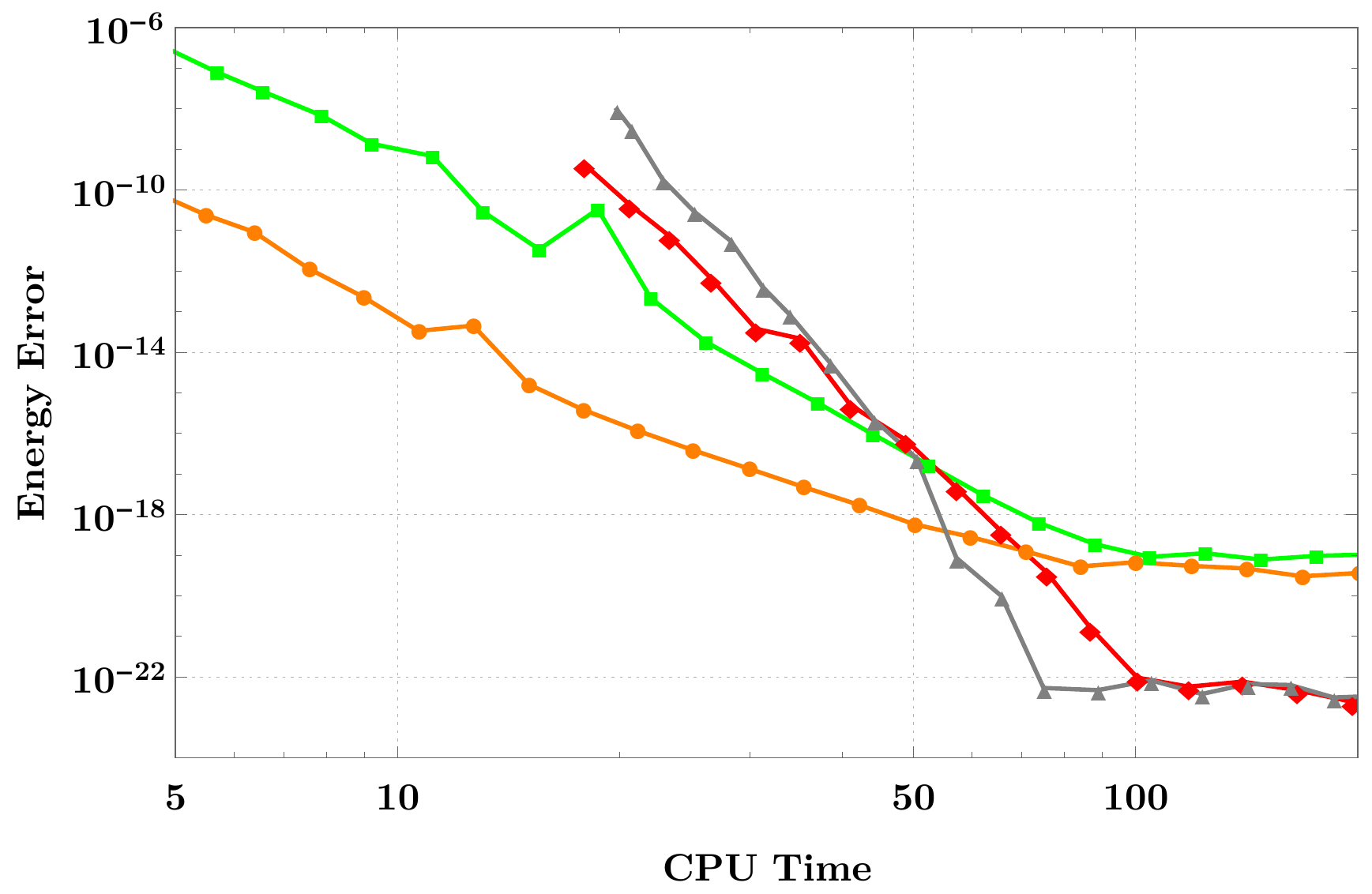}}
&
\subfloat[Parallel executions]
{\includegraphics[width=.45\textwidth]{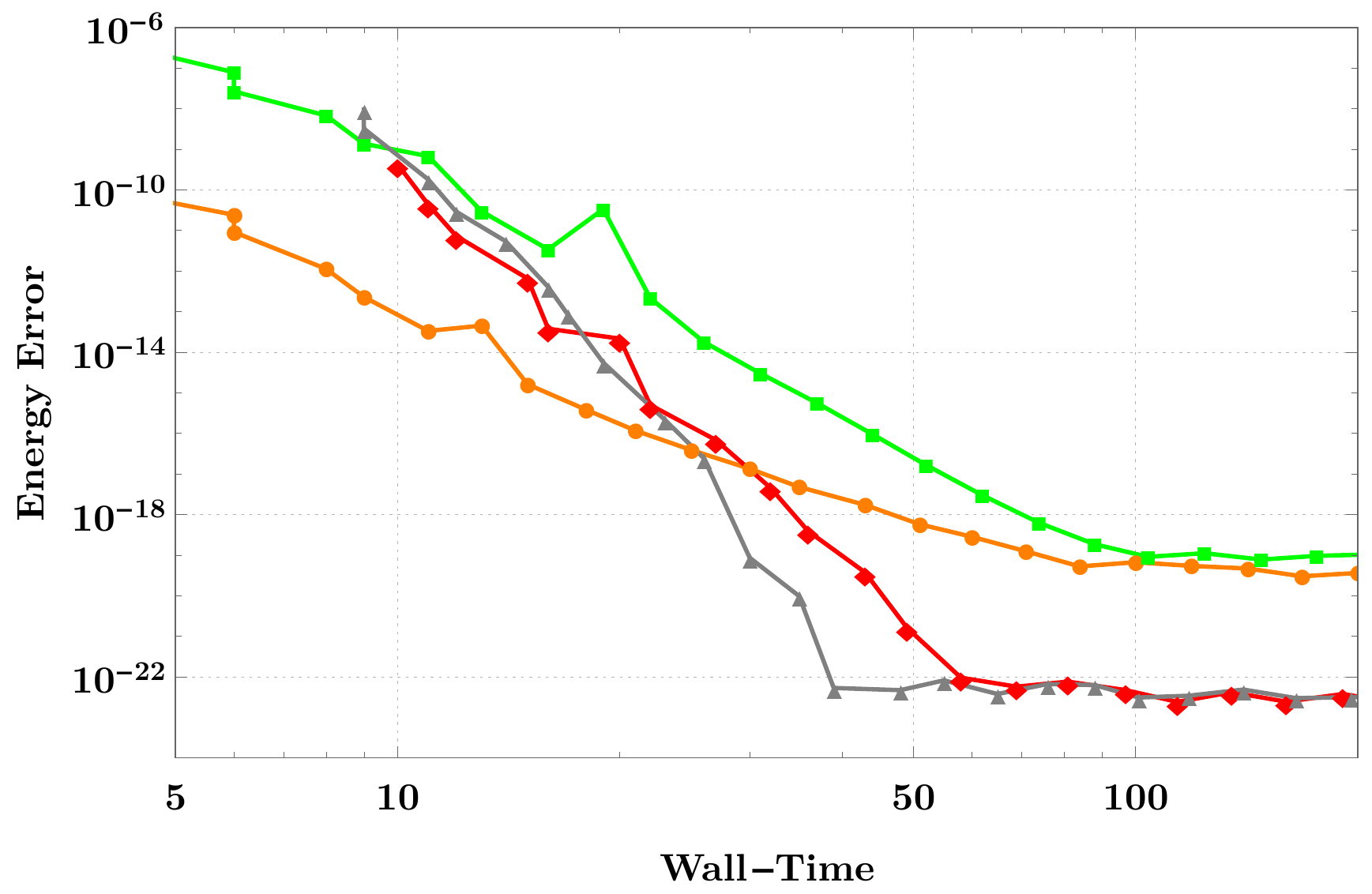}}
\end{tabular}
\caption{\small Efficiency diagrams (in double logarithmic scale) of (a) sequential executions cpu time  and (b) parallel executions wall time using threads=$4$ versus the maximum error in energy. We compare FCIRK methods with $s=8$ stages (red) and $s=16$ stages (gray) implemented with mixed precision (long double - quadruple) and explicit integrators \emph{ABAH1064} (orange) and \emph{CO1035} (green) implemented in long double ($80$-bit) precision arithmetic}
\label{fig:esp4b}
\end{figure}

The mixed precision strategy implemented for FCIRK methods allows us for additional digits of precision when combining long double arithmetic and quadruple precision arithmetic. 
 Despite the comparatively high cpu time cost of evaluating the required flows of the unperturbed system in quadruple precision, the efficiency diagrams of FCIRK shown in  Figures~\ref{fig:esp4a} and \ref{fig:esp4b} are very favorable for the highest precisions.  Indeed,  linear extrapolation of the efficiency diagram of  \emph{ABAH1064}  below its round-off error level shows that the higher order of precision of the FCIRK integrators (which are of order $16$ and $32$ respectively) give them a clear advantage over \emph{ABAH1064} method for the highest precision. In addition, achieving such precisions with  \emph{ABAH1064} would also require a mixed precision implementation (with all the required flows of the unperturbed system obtained in quadruple precision). Moreover, this would  add a considerable penalty in computing times, since several flows have to be evaluated at each integration step of  \emph{ABAH1064} (nine per step, compared to one or two per step of FCIRK, the later using larger, and hence fewer steps, for a given accuracy).

\subsection{Statistical analysis of propagation of round-off errors}

Next we present the result of some numerical experiments designed to asses our implementation of FCIRK integrators with respect to round-off error propagation. We adopt  (as in~\cite{Hairer2008,antonana2017}) a statistical approach. 

The components of the angular momentum 
are quadratic first integrals of both the unperturbed system and the original system under consideration. Recall that FCIRK integrators based on symplectic IRK schemes conserve all quadratic first integrals.  Let us denote as $I(Q,V)$ the Euclidean norm of the angular momentum for given values of the state variables $Q,V \in \mathbb{R}^{3N}$. Hence, $I(Q,V)$ would be exactly conserved by such FCIRK integrators if exact arithmetic were used and the implicit equations required to solve at each step were solved exactly. In practice, this is of course not possible, and we will loosely refer to the actual errors due to the finite precision arithmetic and the finite number of fixed point iterations as round-off errors.

The energy $H(Q,V)$ is also a conserved quantity of the system. Although it will not be exactly conserved by a symplectic FCIRK scheme, the  error in energy will be dominated by the influence of round-off errors  for sufficiently fine time discretization.

We have considered  $P=1000$ perturbed initial values by randomly perturbing each component of the initial values with a relative error of size $\mathcal{O}(10^{-6})$. For each of these initial values, we have integrated the corresponding problem for an integration interval of $10^6$ days, and saved the numerical approximations of $(Q(t),V(t))$ for every $m$ steps of length $h$ at intermediate times $t_{m k} = k m h, \ k=0,1,2,\dots$ obtained by the application of the $32$th order FCIRK scheme in mixed precision (long double - quadruple) arithmetic. 

The local error in angular momentum and energy due to round-off errors, is  "expected" to behave, for a good implementation free from biased errors,  like an independent random variable.  Then, provided that the numerical results are sampled every $m$ steps, 
with a large enough sampling frequency $m$, the difference in energy (resp. norm of angular momentum) will behave as an independent variable with an approximately Gaussian distribution with mean $\mu$ (ideally $\mu=0$) and standard deviation $\sigma$.  So that the accumulated difference in energy (resp. norm of angular momentum)
at the sampled times $t_{m k} = k m h$ would  behave like a Gaussian random walk with standard deviation  $k^{\frac12}\sigma = (t_{m k}/(m h))^{1/2} \sigma$.  This is sometimes referred to as Brouwer's law in the scientific literature~\cite{Grazier2005}, from the original work on the accumulation of round-off errors in the numerical integration of Kepler's problem done by Brouwer in~\cite{Brouwer1937}. For numerical integrators implemented with compensated summation, $\sigma$ is expected to be approximately proportional to $h$. 

Hence, we arrive to the conclusion that ideally, the  relative errors 
\begin{equation*}
I(Q(t),V(t)) /I(Q(0),V(0))-1 
\end{equation*}
in the Euclidean norm $I(Q,V)$  of the angular momentum, and the relative errors in energy
\begin{equation*}
H(Q(t),V(t)) /H(Q(0),V(0))-1 
\end{equation*}
for integrations with small enough $h$ is approximately proportional to $\sqrt{ht}$. 

\begin{figure} [h!]
\centerline{\includegraphics [width=8cm, height=5cm] {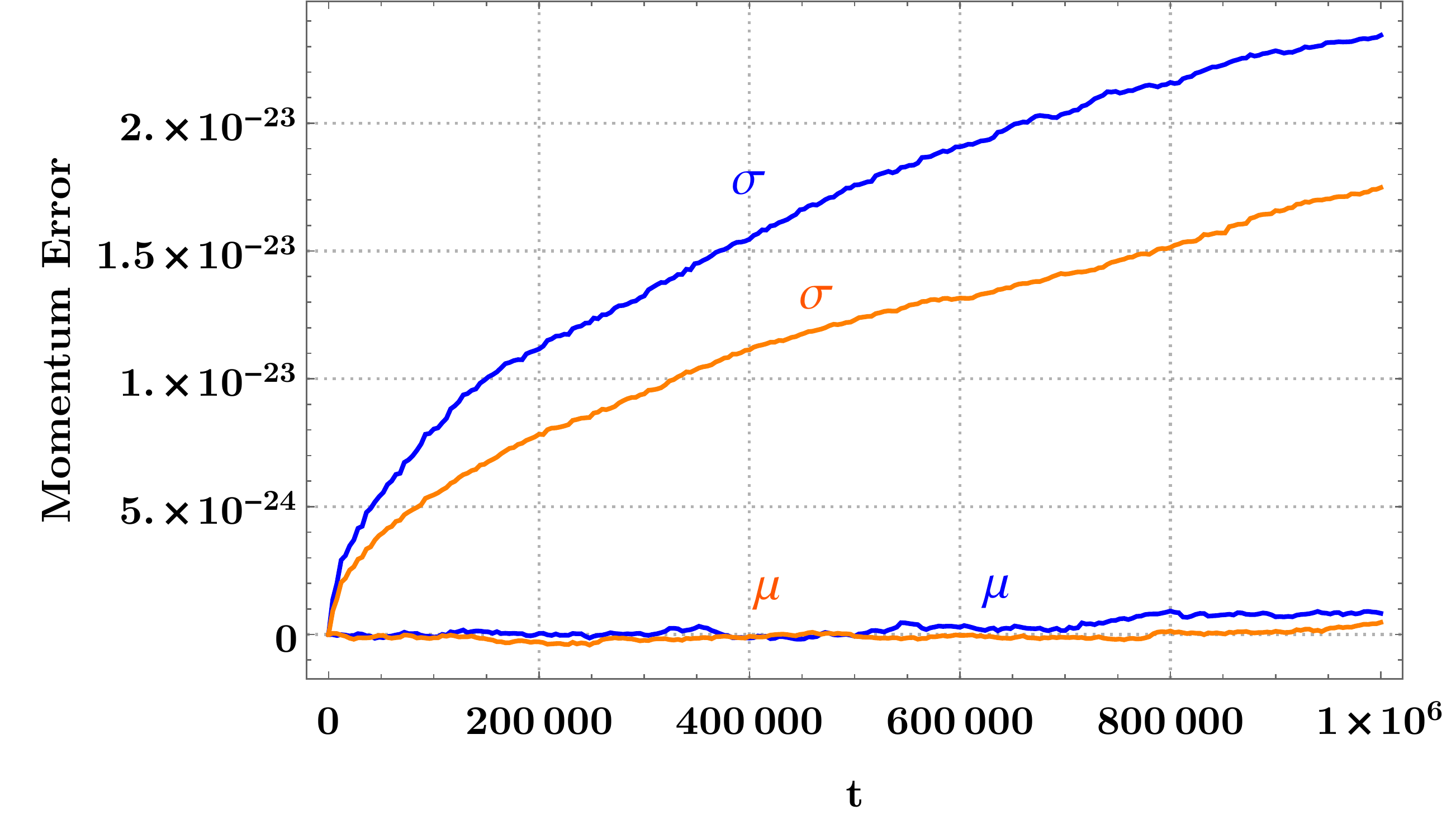}}
\caption{\small Evolution of mean ($\mu$) and standard deviation ($\sigma$) of the relative errors in the Euclidean norm $I(Q,V)$ of the angular momentum: one computation using time-step $h=10$ days (orange) and second one using time-step $h=20$ days (blue)}
\label{fig:esp5}
\end{figure}

In order to check that for the angular momentum, we have made the integrations with two different time-step lengths, one of 10 days, and the other one of 20 days.
In Figure~\ref{fig:esp5} we plot the evolution of the mean and standard deviation of the relative errors in $I$  corresponding to each step length. We observe that in both cases  ($h=10$ and $h=20$ respectively), the standard deviation grows roughly speaking proportionally to $\sqrt(t)$ as time $t$ increases. In addition, it seems that the standard deviation is also proportional to $\sqrt{h}$.  On the other hand, the evolution of the mean errors shows no clear drift as time increases, and seems to be of similar size for both time discretization, and substantially smaller than the standard deviations.

\begin{figure} [h!]
\centerline{\includegraphics [width=8cm, height=5cm] {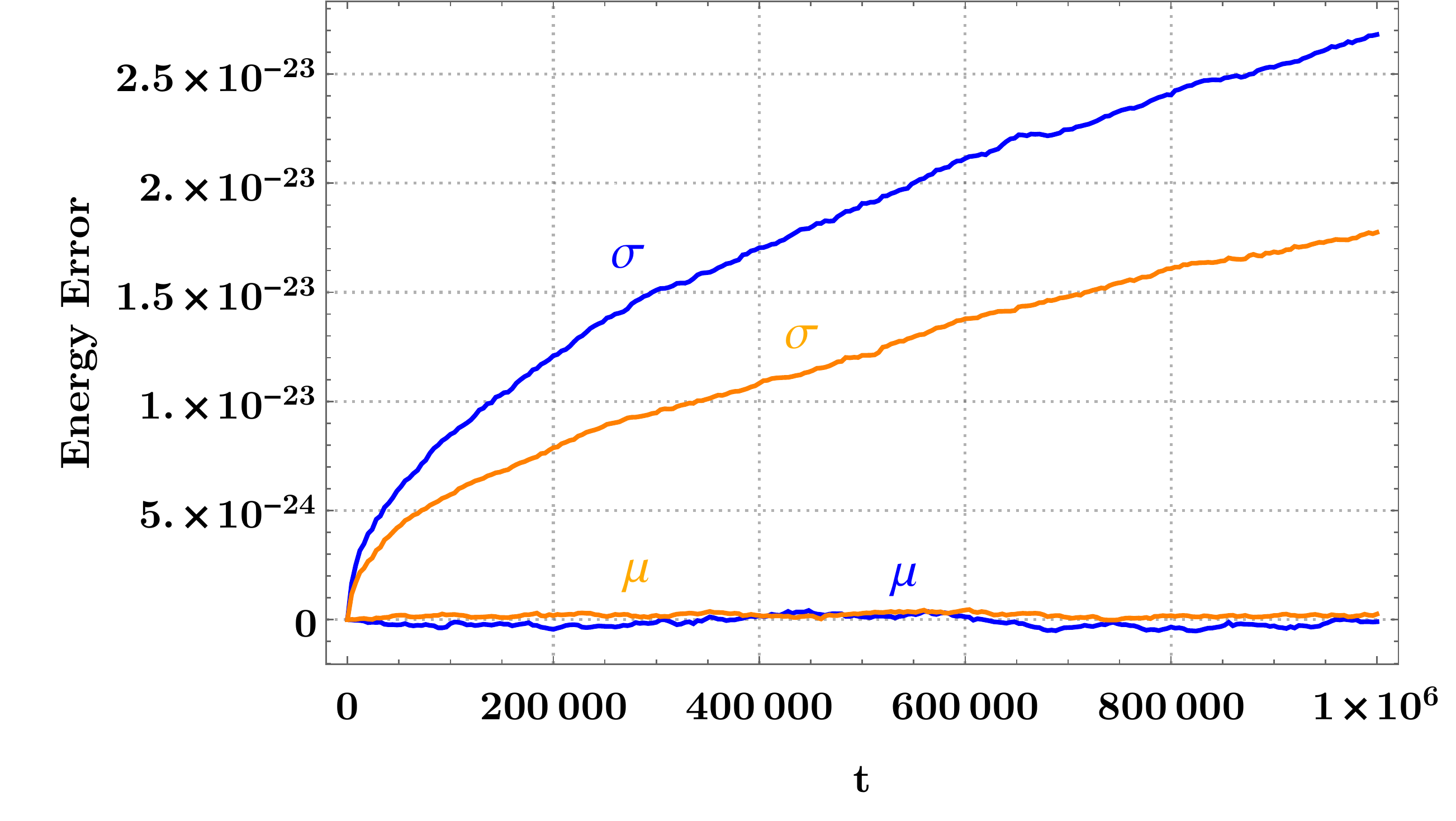}}
\caption{\small Evolution of mean ($\mu$) and standard deviation ($\sigma$) of the relative errors in in energy: one computation using time-step $h=2$ days (orange) and second one using time-step $h=4$ days (blue)}
\label{fig:esp6}
\end{figure}

For the error in energy, we considered finer time discretization in order that the truncation errors be dominated by round-off errors, namely, $h=4$ and $h=2$. 
In Figure~\ref{fig:esp6} we plot the evolution of the mean and standard deviation of the relative errors in energy  corresponding to each step length. Again, the results seem to support the robustness of our implementation with respect to the influence of round-off errors in energy.

\section{Concluding Remarks}
\label{s:cr}

From our numerical experiments, we can conclude  that family of FCIRK schemes of arbitrarily high order has a clear advantage  over splitting methods  for high precision computations (with quadruple precision and arbitrary precision arithmetic).  Even for relatively lower accuracy computations (made in mixed precision arithmetic as in the numerical experiments presented in Section~\ref{s:ne}),  FCIRK methods can potentially outperform splitting integrators, if the underlying IRK scheme is implemented in a parallel computing environment (the $s$ evaluation of the right-hand side of the ODE at each iteration ) and the iteration process is optimized using additional techniques. In particular, one can compute most iterations with a simplified model performed in a lower precision arithmetic, and  at the end compute a few iterations with the full model  in high precision arithmetic to adjust the last digits of the solutions \cite{Beylkin2014}. This seems to be particularly so for more realistic (and hence more complex) models (with more bodies, relativistic effects and other phenomenon taking into account...) where the relative overhead of parallelization  and application of more sophisticated iteration strategies will be less important. 

In addition, splitting methods only allow us to construct symplectic explicit integrators if the perturbation is such that its flow can be explicitly computed, or at least, admits an explicit symplectic approximate discrete flow. Unlike splitting methods,  FCIRK schemes admit any kind of perturbing terms.
On the other hand, although in principle splitting methods admit a mixed precision implementation similar to that mentioned above for FCIRK schemes (the flows of the unperturbed system in higher precision and the evaluations of the perturbation in a lower precision arithmetic), it is computationally more expensive, since several computations of flows of the unperturbed system are required at each step of high order splitting integrators.

\begin{acknowledgements}
M. Anto\~nana, J. Makazaga, and A. Murua have received funding from the Project of the Spanish Ministry of Economy and Competitiveness with reference MTM2016-76329-R (AEI/FEDER, EU), from the
project MTM2013-46553-C3-2-P from Spanish Ministry of Economy and Trade, 
and as part of the Consolidated Research Group IT649-13 by the Basque Government.
\end{acknowledgements}

\bibliographystyle{spmpsci}      

\end{document}